\documentstyle[12pt]{article}

\makeatletter

\def\LRA#1#2{\@tempdimb=\c@enumiv\@tempdima%
   \vcenter{\offinterlineskip\halign{##\cr%
   \hfil${\scriptstyle{#1}}$\hfil\crcr%
   \hbox to \@tempdimb{\rightarrowfill}\cr%
   \noalign{\kern-1ex}%
   \hbox to \@tempdimb{\leftarrowfill}\cr%
   \hfil${\scriptstyle{#2}}$\hfil\crcr}}}
\def\RA#1{\@tempdimb=\c@enumiv\@tempdima\vbox{\offinterlineskip%
   \halign{##\cr\hfil${\scriptstyle {#1}}$\hfil\crcr%
   \hbox to \@tempdimb{\rightarrowfill}\cr}}}
\def\LA#1{\@tempdimb=\c@enumiv\@tempdima\vbox{\offinterlineskip%
   \halign{##\cr\hfil${\scriptstyle {#1}}$\hfil\crcr%
   \hbox to \@tempdimb{\leftarrowfill}\cr}}}

\def\diag{\leavevmode\bgroup\setcounter{enumiv}{1}%
   \unitlength1em \@tempdima3em \def\\{\crcr&}\vbox\bgroup%
   \def\multicolumn##1##2{\multispan##1\setcounter{enumiv}{##1}%
   \hfil{##2}\hfil\setcounter{enumiv}{1}}
   \offinterlineskip\halign\bgroup\vrule height.8em depth.7em %
   width0pt##&&\hfil${\displaystyle{##}}$\hfil\cr&}
\def\enddiag{\crcr\egroup\egroup\egroup}
\makeatother
\if@twoside \oddsidemargin -20pt \evensidemargin 20pt \marginparwidth 85pt
\else \oddsidemargin 0pt \evensidemargin 0pt
\fi
\advance\textheight by \topskip
\font\symbolfont=msbm10
\font\teneu=eufm10
\font\egteu=eufm8
\def\dn#1{\mathchoice{\hbox{\teneu #1}}{\hbox{\teneu #1}}%
   {\hbox{\egteu #1}}{\hbox{\egteu #1}}}

\def\eee{\hbox{\symbolfont E}}
\def\fff{\hbox{\symbolfont F}}

\def\nnn{\hbox{\symbolfont N}}
\def\zzz{\hbox{\symbolfont Z}}
\def\ppp{\hbox{\symbolfont P}}
\def\qqq{\hbox{\symbolfont Q}}
\def\rrr{\hbox{\symbolfont R}}

\def\ttt{\hbox{\symbolfont T}}
\def\iii{\hbox{\symbolfont I}}

\def\cc{\mathop{\cal C}\nolimits}
\def\dd{\mathop{\cal D}\nolimits}

\def\pp{\mathop{\cal P}\nolimits}
\def\qq{\mathop{\cal Q}\nolimits}
\def\ff{\mathop{\cal F}\nolimits}
\def\tt{\mathop{\cal T}\nolimits}

\def\jac{\mathop{\cal J}\nolimits}
\def\pd{\mathop{\rm pd}\nolimits}

\def\id{\mathop{\rm id}\nolimits}
\def\gl{\mathop{\rm gl.dim.}\nolimits}

\def\spep#1{\mathop{{}^{\bullet}\strut\kern-.1em{#1}}\nolimits}
\def\gd{\mathop{\rm Gr}\nolimits}
\def\ext{\mathop{\rm Ext}\nolimits}
\def\hom{\mathop{\rm Hom}\nolimits}
\def\endm{\mathop{\rm End}\nolimits}

\def\add{\mathop{\rm add}\nolimits}

\def\soc{\mathop{\rm soc}\nolimits}
\def\top{\mathop{\rm top}\nolimits}

\def\pr{\mathop{\rm pr}\nolimits}
\def\inj{\mathop{\rm in}\nolimits}
\def\ssim{\mathop{\rm sim}\nolimits}
\def\ind{\mathop{\rm ind}\nolimits}

\def\ann{\mathop{\rm ann}\nolimits}
\def\tr{\mathop{\rm Tr}\nolimits}

\def\Cok{\mathop{\rm Cok}\nolimits}

\def\Im{\mathop{\rm Im}\nolimits}

\def\supp{\mathop{\rm supp}\nolimits}

\def\nak{\mathop{\rm n}\nolimits}

\def\mod{\mathop{\rm mod}\nolimits}

\def\extf#1#2{\mathop{\eee_{#1}#2}}

\def\add{\mathop{\rm add}\nolimits}

\def\gd{\mathop{\rm Gr}\nolimits}

\def\cc{\mathop{\cal C}\nolimits}

\def\Cok{\mathop{\rm Cok}\nolimits}

\def\pd{\mathop{\rm pd}\nolimits}
\def\fd{\mathop{\rm fd}\nolimits}
\def\id{\mathop{\rm id}\nolimits}
\def\gl{\mathop{\rm gl.dim}\nolimits}
\def\domdim{\mathop{\rm dom.dim}\nolimits}

\def\ext{\mathop{\rm Ext}\nolimits}

\def\pp{\mathop{\ind^+_1}\nolimits}
\def\ii{\mathop{\ind^-_1}\nolimits}
\def\qq{\mathop{\cal Q}\nolimits}

\def\tr{\mathop{\rm Tr}\nolimits}
\def\length{\mathop{\rm length}\nolimits}

\def\grade#1{\mathop{{\rm grade}\ #1}}
\def\sgrade#1{\mathop{{\rm s.grade}\ #1}}
\def\rgrade#1{\mathop{{\rm r.grade}\ #1}}

\def\AR{\hbox{\symbolfont A}}

\def\quo{\mathop{\rm Fac}\nolimits}
\def\XA{1}
\def\XAA{1.1}
\def\XAB{1.2}
\def\XABA{1.2.1}
\def\XAC{1.3}
\def\XACA{1.3.1}
\def\XACB{1.3.2}
\def\XACC{1.3.3}
\def\XAD{1.4}
\def\XADA{1.4.1}
\def\XADB{1.4.2}
\def\XAE{1.5}
\def\XB{2}
\def\XBA{2.1}
\def\XBAA{2.1.1}
\def\XBB{2.2}
\def\XBBA{2.2.1}
\def\XBBB{2.2.2}
\def\XBBC{2.2.3}
\def\XBC{2.3}
\def\XBD{2.4}
\def\XBE{2.5}
\def\XBEA{2.5.1}
\def\XC{3}
\def\XCA{3.1}
\def\XCAA{3.1.1}
\def\XCAB{3.1.2}
\def\XCB{3.2}
\def\XCC{3.3}
\def\XCCA{3.3.1}
\def\XCD{3.4}
\def\XCDA{3.4.1}
\def\XCDB{3.4.2}
\def\XCE{3.5}
\def\XCF{3.6}
\def\XCFA{3.6.1}
\def\XCFB{3.6.2}
\def\XCFC{3.6.3}
\def\XCFD{3.6.4}
\def\XD{4}
\def\XDA{4.1}
\def\XDB{4.2}
\def\XDBA{4.2.1}
\def\XDBB{4.2.2}
\def\XDBC{4.2.3}
\def\XDC{4.3}
\def\XDD{4.4}
\def\XDDA{4.4.1}
\def\XDDB{4.4.2}
\def\XDE{4.5}
\def\XE{5}
\def\XEA{5.1}
\def\XEAA{5.1.1}
\def\XEB{5.2}
\def\XEC{5.3}
\def\XECA{5.3.1}
\def\XED{5.4}
\def\XEE{5.5}
\def\XEEA{5.5.1}
\def\Q#1{\begin{picture}(0,0)\put(-1,-2){#1}\end{picture}}
\def\QUR{\begin{picture}(4,0)\put(0,-5){\vector(1,1){10}}\end{picture}}
\def\QDR{\begin{picture}(4,0)\put(0,5){\vector(1,-1){10}}\end{picture}}
\def\QR{\begin{picture}(4,0)\put(0,0){\vector(1,0){10}}\end{picture}}
\begin{document}
\begin{center}
\vspace*{1cm}{\large The relationship between homological properties and representation theoretic realization of artin algebras}
\vskip1em{\sc Osamu Iyama}
\end{center}

{\footnotesize
{\sc Abstract. }We will study the relationship of quite different object in the theory of artin algebras, namely Auslander-regular rings of global dimension two, torsion theories, $\tau$-categories and almost abelian categories. We will apply our results to characterization problems of Auslander-Reiten quivers.}

\vskip1em{\bf 0.1 }
There exists a bijection between equivalence classes of Krull-Schmidt categories $\cc$ with additive generators $M$ and Morita-equivalence classes of semiperfect rings $\Gamma$, which is given by $\cc\mapsto\cc(M,M)$ and the converse is given by $\Gamma\mapsto\pr\Gamma$ for the category $\pr\Gamma$ of finitely generated projective $\Gamma$-modules. Although this bijection itself is rather formal, it will be very fruitful to study the relationship of (A)--(D) below. The object of this paper is to study it under the assumption that $\Gamma$ is an artin algebra.

\strut\kern1em
(A) Homological properties for $\Gamma$.

\strut\kern1em
(B) Representation theoretic realization of $\cc$.

\strut\kern1em
(C) Categorical properties for $\cc$.

\strut\kern1em
(D) Combinatorial properties for the AR quiver $\AR(\cc)$.

For (A), we will study a property of the selfinjective resolution of $\Gamma$ which is called the {\it $(l,n)$-conditions} (\S\XAA) and generalizes both the Auslander conditions [Bj] and the dominant dimension [T]. For (B), we will study the existence of an equivalence between $\cc$ and a torsionfree class of $\mod\Lambda$ over an artin algebra $\Lambda$ (\S\XAB,\S\XBB), where such a subcategory is very popular in the representation theory of artin algebras [Ha][As]. For (C), we will treat a class of additive categories which are called {\it $\tau$-categories} (\S\XAC) and introduced in [I3]. $\tau$-categories are additive categories with generalized almost split sequences, and our motivation and definition were rather different from the work of Auslander and Smalo in [AS] since we aimed to treat categories which can be far from abelian, for example mesh categories of translation quivers (\S\XACB(3)). Nevertheless our result \S\XBA\ asserts that some $\tau$-categories are realized as torsionfree classes over artin algebras, and they form almost abelian categories (\S\XAE). For (D), we will study a combinatorial invariant $\AR(\cc)$ of a $\tau$-category $\cc$ called the {\it AR} (=Auslander-Reiten) {\it quiver} (\S\XDA). Some results in this paper were already announced in [I6;7.4] without proof.

\vskip1em{\bf 0.2 Background }
In [A], Auslander obtained a quite remarkable theorem which asserts that there exists a bijection between Morita-equivalence classes of representation-finite artin algebras $\Lambda$ and those of {\it Auslander algebras} $\Gamma$, which is an artin algebra with $\gl\Gamma\le2$ and $\domdim\Gamma\ge2$. Then $\mod\Lambda$ is equivalent to $\pr\Gamma$, and this correspondence gives a prototype of our study. It relates a representation theoretic realization (B) of the category $\cc=\pr\Gamma$ to a homological property (A) of $\Gamma$, and it will be suggestive that Auslander algebras form a special class of Auslander-regular rings $\Gamma$ with $\gl\Gamma\le2$. In [FGR][Bj][AR2][C] and so on, Auslander-regular rings, more generally Auslander-Gorenstein rings, are studied as a non-commutative analogy of commutative Gorenstein rings motivated by the classical results of Bass [B].

Later, Auslander and Reiten [AR1] obtained the existence theorem of almost split sequences, which is one of the most important theorems in the representation theory of algebras [ARS]. This theorem gives a categorical property (C) of $\mod\Lambda$ over an artin algebra $\Lambda$, which means that `$\mod\Lambda$ forms a $\tau$-category' in our context. Although this theorem is a great achivement of general theory of homological algebra, it can be proved easily for the special case when $\Lambda$ is representation-finite, and this observation seems to lead them to the general existence theorem. Moreover, the theorem for representation-finite case has its own importance. In terms of the corresponding Auslander algebra $\Gamma$, it means that the functor $\ext^2_\Gamma(\ ,\Gamma):\mod\Gamma\leftrightarrow\mod\Gamma^{op}$ gives a duality between simple $\Gamma$-modules $L$ with $\pd L=2$ and that of $\Gamma^{op}$. We can naturally generalize this duality to arbitrary Auslander-Gorenstein rings [I9] (see \S\XCFC\ below).

By applying the Auslander's correspondence theorem and extending some aspects in [BG], Igusa-Todorov and Brenner gave (distinct) characterizations of AR quivers of representation-finite artin algebras [IT3][Br]. These are nothing but combinatorial properties (D) of $\AR(\cc)$. Recently, inspired by the work of Igusa-Todorov, the author introduced $\tau$-categories and successfully applied them to characterize AR quivers of representation-finite orders [I3,4,5]. We shall see that $\tau$-categories give a powerful tool for our problem.

\vskip1em{\bf 0.3 Our results }
Our first theorem \S\XBA\ gives the relationship of (A)--(D) in \S0.1 for more general classes of algebras than those in \S0.2, namely the conditions below are equivalent for an artin algebra $\Gamma$ and $\cc=\pr\Gamma$.

\strut\kern1em
(A) $\Gamma$ satisfies $\gl\Gamma\le2$ and the two-sided $(2,2)$-condition (\S\XAA).

\strut\kern1em
(B) $\cc$ is equivalent to a faithful torsionfree class over an artin algebra (\S\XAB).

\strut\kern1em
(C) $\cc$ is a strict $\tau$-category (\S\XAC).

\strut\kern1em
(D) $\cc$ is a $\tau$-category with a right additive function on $\AR(\cc)$ (\S\XACA).

Next we will study Auslander-regular rings $\Gamma$ with $\gl\Gamma\le2$, which forms a special class of algebras in above (A). Our second theorem \S\XCA\ gives the corresponding objects in (B)--(D) to such $\Gamma$, namely hereditary torsionfree classes (\S\XAB) correspond for (B), strict $\tau$-categories with `Nakayama pairs' (\S\XAD) correspond for (C), and $\tau$-categories with additive functions on $\AR(\cc)$ (\S\XACA) correspond for (D). The concept of Nakayama pairs was introduced in [I4] to characterize AR quivers of representation-finite orders, and they were essentially used also in Igusa-Todorov's theorem (\S0.2).
The concept of additve functions often appeared in the representation theory (see \S\XD), for example Brenner's theorem (\S0.2) and Rump's recent characterization of AR quivers of representation-finite orders [R5]. Moreover, we will study two special classes of Auslander-regular rings $\Gamma$ with $\gl\Gamma\le2$, and give the corresponding objects in (B)--(D) again. One is Auslander algebras (\S\XCC) which gave prototype of our study, and we will prove very clearly Auslander's correspondence theorem, Igusa-Todorov's theorem and Brenner's theorem explained in \S0.2.
Another is diagonal Auslander-regular rings $\Gamma$ with $\gl\Gamma\le2$ (\S\XCD), which are closely related to the category $\mod_{sp}\Lambda$ of modules with projective socles [S], and we will generalize a result of Ringel-Vossieck [RV]. The connection of several known results will be understood clearly in our viewpoint of this paper. As we shall see in examples in \S\XDE, our characterizations of AR quivers \S\XDD\ are very simple and can be checked easily. 

We will discuss properties of $\Gamma$ with $\gl\Gamma\le2$ and the two-sided $(2,2)$-condition, namely symmetry in \S\XCFA, duality in \S\XCFC\ and quasi-Koszul property [GM] of Green-Martinez in \S\XBE.
In the final section \S\XE, we will study the rejection theory of faithful torsionfree classes over artin algebras. The rejection theory of $\tau$-categories was studied in [I1,4] as a wide generalization of DK (=Drozd-Kirichenko) Rejection Lemma which was fundamental in the theory of Bass orders ([DKR][Ro][HN]). They played a crucial role to characterize AR quivers of representation-finite orders [I5], and recently they were applied to prove Solomon's second conjecture on zeta functions of orders [I7] and finiteness of representation dimension of artin algebras [I8].

All artin algebras $\Gamma$ in (A) studied in this paper satisfies $\gl\Gamma\le2$. It will be very interesting to generalize our results to artin algebras $\Gamma$ with $\gl\Gamma\ge3$.

Finally, notice that W. Rump's recent work [R1,2,3,4,5] on almost abelian categories (\S\XAE) and $\tau$-categories has strong relationship with our study.

\vskip1em{\bf\XA\ Preliminaries }

In this paper, any module is assumed to be a left module.
For a ring $\Lambda$, we denote by $\mod\Lambda$ the category of finitely generated $\Lambda$-modules, by $\pr\Lambda$ (resp. $\ssim\Lambda$) the category of finitely generated projective (resp. simple) $\Lambda$-modules, by $J_\Lambda$ the Jacobson radical of $\Lambda$, by $\widehat{(\ )}$ the functor $\hom_\Lambda(\ ,\Lambda):\mod\Lambda\leftrightarrow\mod\Lambda^{op}$, by $0\rightarrow X\stackrel{}{\rightarrow}I^0_\Lambda(X)\stackrel{}{\rightarrow}I^1_\Lambda(X)\stackrel{}{\rightarrow}\cdots$ a minimal injective resolution of a $\Lambda$-module $X$, and by $\pd X$ (resp. $\fd X$, $\id X$) the projective (resp. flat, injective) dimension of a $\Lambda$-module $X$.
We denote by $\underline{\mod}\Lambda$ the stable category, by $\Omega:\underline{\mod}\Lambda\rightarrow\underline{\mod}\Lambda$ the syzygy functor, and by $\tr:\underline{\mod}\Lambda\rightarrow\underline{\mod}\Lambda^{op}$ the transpose functor [AB].
When $\Lambda$ is an artin algebra over $R$, we denote by $(\ )^*$ the duality $\hom_R(\ ,I^0_R(R/J_R)):\mod\Lambda\leftrightarrow\mod\Lambda^{op}$, and by $\inj\Lambda$ the category of finitely generated injective $\Lambda$-modules.
For an additive category $\cc$ and $X\in\cc$, we denote by $\add X$ the full subcategory of $\cc$ consisting of direct summands of $X^n$ ($n>0$). We call $X$ an {\it additive generator} of $\cc$ if $\cc=\add X$ holds.

\vskip1em{\bf\XAA\ }
Let $\Gamma$ be a noetherian ring.
We denote by $\grade{L}:=\inf\{ i\ge0\ |\ \ext_\Gamma^i(L,\Gamma)\neq0\}$ (resp. $\sgrade{L}:=\inf\{ \grade{M}\ |\ M\subseteq L\}$, $\rgrade{L}:=\inf\{ i>0\ |\ \ext_\Gamma^i(L,\Gamma)\neq0\}$) the {\it grade} (resp. {\it strong grade}, {\it reduced grade}) of $L\in\mod\Gamma$. For any $n\ge0$, the full subcategory $\{ L\ |\ \sgrade{L}\ge n\}$ of $\mod\Gamma$ is abelian and closed under subfactor modules and extensions.
For $l$, $n>0$, we say that $\Gamma$ satisfies the {\it $(l,n)$-condition} if the following equivalent conditions are satisfied [I5;6.1].

\strut\kern1em(i) $\fd I^i_\Gamma(\Gamma)<l$ holds for any $i$ ($0\le i<n$).

\strut\kern1em(ii) $\sgrade{\ext_\Gamma^l(L,\Gamma)}\ge n$ holds for any $L\in\mod\Gamma^{op}$.

This equivalence much simplifies the equivalence of (a) and (c) in the famous theorem [FGR;3.7] and that of (b) and (d) in [AR3;0.1]. For an artin algebra $\Gamma$, the $(l,n)$-condition is equivalent to that $\sgrade{\ext_\Gamma^l(L,\Gamma)}\ge n$ holds for any {\it simple} $\Gamma^{op}$-module $L$ [AR2]. We say that $\Gamma$ satisfies the {\it $(l,n)^{op}$-condition} if $\Gamma^{op}$ satisfies the $(l,n)$-condition.

Put $\domdim\Gamma:=\inf\{ i\ge0\ |\ \fd I^i_\Gamma(\Gamma)\neq0\}$ [T], which is the maximal number $n$ such that $\Gamma$ satisfies the $(1,n)$-condition. We call $\Gamma$ {\it $n$-Gorenstein} if $\fd I^i_\Gamma(\Gamma)\le i$ holds for any $i$ ($0\le i<n$) [FGR], or equivalently, $\Gamma$ satisfies the $(l,l)$-condition for any $l$ ($0<l\le n$). We call $\Gamma$ {\it Auslander-regular} (resp. {\it Auslander-Gorenstein}) if $\gl\Gamma<\infty$ (resp. $\id{}_\Gamma\Gamma<\infty$ and $\id\Gamma_\Gamma<\infty$) and $\Gamma$ is $n$-Gorenstein for any $n$ [C]. It is well-known that $\domdim\Gamma=\domdim\Gamma^{op}$ holds [H2], and $\Gamma$ is $n$-Gorenstein if and only if so is $\Gamma^{op}$ [FGR;3.7]. These left-right symmetries were generalized to the $(l,n)$-condition in [I9] (see \XCFA\ below) though the $(l,n)$-condition itself is not left-right symmetric \XBAA(2).

\vskip1em{\bf\XAB\ }
Let $\Lambda$ be an artin algebra. We call a full subcategory $\cc$ of $\mod\Lambda$ a {\it torsionfree} (resp. {\it torsion}) {\it class} if it is closed under submodules (resp. factor modules) and extensions [As]. For a collection ${\bf S}$ of $\Lambda$-modules, define full subcategories of $\mod\Lambda$ by ${\bf S}^\perp:=\{ X\ |\ \hom_\Lambda(Y,X)=0$ for any $Y\in{\bf S}\}$ and ${}^\perp{\bf S}:=\{ X\ |\ \hom_\Lambda(X,Y)=0$ for any $Y\in{\bf S}\}$. We call a pair $(\tt,\ff)$ of full subcategories of $\mod\Lambda$ a {\it torsion theory} on $\mod\Lambda$ if it satisfies the following equivalent conditions:

\strut\kern1em(i) $\ff=\tt^\perp$ and $\tt={}^\perp\ff$.

\strut\kern1em(ii) $\ff$ is a torsionfree class and $\tt={}^\perp\ff$.

\strut\kern1em(iii) $\tt$ is a torsion class and $\ff=\tt^\perp$.

\strut\kern1em(iv) The inclusion $\ff\rightarrow\mod\Lambda$ has a left adjoint $\fff:\mod\Lambda\rightarrow\ff$ with a unit $\alpha$ and the inclusion $\tt\rightarrow\mod\Lambda$ has a right adjoint $\ttt:\mod\Lambda\rightarrow\tt$ with a counit $\beta$ such that $0\rightarrow\ttt\stackrel{\beta}{\rightarrow}1_{\mod\Lambda}\stackrel{\alpha}{\rightarrow}\fff\rightarrow0$ is exact.

We call a torsion theory $(\tt,\ff)$ (resp. torsion class $\tt$, torsionfree class $\ff$) {\it faithful} if $\Lambda\in\ff$, {\it cofaithful} if $\Lambda^*\in\tt$, {\it hereditary} if $\tt$ is closed under submodules, and {\it cohereditary} if $\ff$ is closed under factor modules.
If $(\tt,\ff)$ is a hereditary torsion theory, then $\ff={\bf S}^{\perp}$ holds for the set ${\bf S}$ of simple $\Lambda$-modules in $\tt$. The facts below show that the faithfulness is fundamental for torsion theories.

\vskip1em{\bf\XABA\ Proposition }{\it
Let $(\tt,\ff)$ and $(\tt_i,\ff_i)$ ($i=1,2$) be torsion theories on $\mod\Lambda$.

(1) There exists a factor algebra $\Gamma$ of $\Lambda$ such that $\ff\subseteq\mod\Gamma$ and $(\tt\cap\mod\Gamma,\ff)$ is a faithful torsion theory on $\mod\Gamma$.

(2) $\ff$ is faithful and contravariantly finite if and only if there exists a cotilting $\Lambda$-module $U\in\mod\Lambda$ with $\id U\le 1$ such that $\tt={}^\perp U$.

(3) Assume that $(\tt_i,\ff_i)$ is faithful for $i=1,2$. Then any equivalence $\iii:\ff_1\rightarrow\ff_2$ extends uniquely to an equivalence $\mod\Lambda_1\rightarrow\mod\Lambda_2$.}

\vskip1em{\sc Proof }
(1) Put $I:=\bigcap_{X\in\ff}\ann_\Lambda X$ and $\Gamma:=\Lambda/I$. Then $\ff\subseteq\mod\Gamma$ holds, and there exists a faithful $\Gamma$-module $Y\in\ff$. Put $E:=\endm_\Gamma(Y)$ and take a surjection $f\in\hom_E(E^n,Y)$. Taking $\hom_E(\ ,Y)$, we obtain an injection $(f\cdot)\in\hom_\Gamma(\endm_E(Y),Y^n)$. Thus $\Gamma\in\ff$ holds by $\Gamma\subseteq\endm_E(Y)$. Obviously $(\tt\cap\mod\Gamma,\ff)$ forms a torsion theory on $\mod\Gamma$.

(2) Well-known (see [H1][AS][As]).

(3) Suppose that $X\in\pr\Lambda_1$ satisfies $\iii X\notin\pr\Lambda_2$. Since $\Lambda_2\in\ff_2$ holds, there exists $f\in\hom_{\Lambda_1}(Y,X)$ such that $\iii f$ is a non-split surjection. Since $f$ also does not split, $Z:=\Im f$ is a proper submodule of $X$. Let $g\in\hom_{\Lambda_1}(Z,X)$ be an injection. Since $\iii g$ is a monomorphism in a torsionfree class $\ff_2$, it is injective in $\mod\Lambda_2$. Thus $\iii f$ factors through a proper submodule $\iii Z$ of $\iii X$, a contradiction. Hence $\iii$ restricts to an equivalence $\pr\Lambda_1\rightarrow\pr\Lambda_2$.
Take a progenerator $X\in\pr\Lambda_1$ such that $\iii X=\Lambda_2$. Since $\hom_{\Lambda_1}(X,\ )\stackrel{\iii}{\stackrel{\sim}{\rightarrow}}\hom_{\Lambda_2}(\iii X,\iii(\ ))=\iii$ holds on $\ff_1$, $\iii$ extends to the equivalence $\hom_{\Lambda_1}(X,\ ):\mod\Lambda_1\rightarrow\mod\Lambda_2$, which is easily shown to be the unique extension of $\iii$.\rule{5pt}{10pt}

\vskip1em{\bf\XAC\ }[I3][R5] Let $\cc$ be a skeletally small additive category. We denote by $\cc(X,Y)$ the set of morphisms from $X$ to $Y$, by $fg$ the composition of $f\in\cc(X,Y)$ and $g\in\cc(Y,Z)$, by $\jac_{\cc}$ the Jacobson radical of $\cc$, and by $\ind\cc$ the set of isoclasses of indecomposable objects in $\cc$. We call $\cc$ {\it Krull-Schmidt} if any object is isomorphic to a finite direct sum of objects whose endomorphism rings are local. 

Let $\cc$ be a Krull-Schmidt category and ${\bf A}:X\stackrel{f}{\rightarrow}Y\stackrel{g}{\rightarrow}Z$ a complex.
We call $f$ a {\it weak-kernel} of $g$ if $\cc(\ ,X)\stackrel{\cdot f}{\rightarrow}\cc(\ ,Y)\stackrel{\cdot g}{\rightarrow}\cc(\ ,Z)$ is exact, and $g$ a {\it weak-cokernel} of $f$ if $\cc(X,\ )\stackrel{f\cdot}{\leftarrow}\cc(Y,\ )\stackrel{g\cdot}{\leftarrow}\cc(Z,\ )$ is exact.
A weak-kernel (resp. weak-cokernel) is called {\it minimal} if it has no direct summand of the form $W\rightarrow0$ (resp. $0\rightarrow W$) with $W\neq0$ as a complex. Clearly a minimal weak-(co)kernel is unique up to isomorphism of complexes if it exists. Now we consider the following conditions for ${\bf A}$.

\strut\kern1em
(i) $f,g\in\jac_{\cc}$, and $0\leftarrow\jac_{\cc}(X,\ )\stackrel{f\cdot}{\leftarrow}\cc(Y,\ )$ and $\cc(\ ,Y)\stackrel{\cdot g}{\rightarrow}\jac_{\cc}(\ ,Z)\rightarrow0$ are exact.

\strut\kern1em
(ii) $f$ is a minimal weak-kernel of $g$.

\strut\kern1em
(iii) $g$ is a minimal weak-cokernel of $f$.

We call ${\bf A}$ a {\it $\tau$-sequence} (resp. {\it right $\tau$-sequence}, {\it left $\tau$-sequence}) if it satisfies (i)(ii)(iii) (resp. (i)(ii), (i)(iii)). We call a right (resp. left) $\tau$-sequence ${\bf A}$ {\it strict} if $f$ is a monomorphism (resp. $g$ is an epimorphism) in $\cc$. They are analogue of almost split sequences in arbitrary Krull-Schmidt categories.

We call $\cc$ a ({\it strict}) {\it $\tau$-category} if any $X\in\cc$ is a right term of some (strict) right $\tau$-sequence and a left term of some (strict) left $\tau$-sequence. Then the right $\tau$-sequence with the right term $X$ and left $\tau$-sequence with the left term $X$ are unique up to isomorphism of complexes, and we sometimes denote them by $(X]=(\tau^+X\stackrel{\nu^+_X}{\rightarrow}\theta^+X\stackrel{\mu^+_X}{\rightarrow}X)$ and $[X)=(X\stackrel{\mu^-_X}{\rightarrow}\theta^-X\stackrel{\nu^-_X}{\rightarrow}\tau^-X)$. We denote by $\pp\cc$ (resp. $\ind^+_0\cc$, $\ii\cc$, $\ind^-_0\cc$) the subset of $\ind\cc$ consisting of $X$ satisfying $\tau^+X=0$ (resp. $\theta^+X=0$, $\tau^-X=0$, $\theta^-X=0$). Up to isomorphism of complexes, $(X]=[\tau^+X)$ and $[Y)=(\tau^-Y]$ hold for any $X\in\ind\cc-\pp\cc$ and $Y\in\ind\cc-\ii\cc$. In particular, $\tau^+$ and $\tau^-$ give mutually inverse bijections between $\ind\cc-\pp\cc$ and $\ind\cc-\ii\cc$ [I3;2.3].

We will use in \XEC\ an important property of $\tau$-categories $\cc$ which asserts that the factor category $\cc/[\cc^\prime]$ forms a $\tau$-category again for any subcategories $\cc^\prime$ of $\cc$, where $[\cc^\prime]$ is the ideal of $\cc$ consisting of morphisms which factor through some object in $\cc^\prime$ [I4;1.4].

\vskip1em{\bf\XACA\ }For a set $Q$, we denote by $\zzz Q$ (resp. $\nnn Q$) the free $\zzz$-module (resp. free abelian monoid) generated by $Q$.
Let $\cc$ be a $\tau$-category. We identify $\nnn\ind\cc$ with the set of isoclasses of objects in $\cc$. We can regard $\theta^+$, $\theta^-$, $\tau^+$ and $\tau^-$ as elements of $\endm_{\zzz}(\zzz\ind\cc)$. Put $\phi^\pm:=1_{\zzz\ind\cc}-\theta^\pm+\tau^\pm\in\endm_{\zzz}(\zzz\ind\cc)$.
Let $l:\ind\cc\rightarrow\nnn_{>0}$ be a map. We extend $l$ uniquely to $l\in\hom_{\zzz}(\zzz\ind\cc,\zzz)$. We call $l$ a {\it right} (resp. {\it left}) {\it additive function} if $l(\phi^+X)\ge0$ (resp. $l(\phi^-X)\ge0$) holds for any $X\in\cc$ and the equality holds for any $X\in\ind\cc-\pp\cc$ (resp. $X\in\ind\cc-\ii\cc$). Then put $l^+:=\{ X\in\pp\cc\ |\ l(\phi^+X)>0\}$ (resp. $l^-:=\{ X\in\ii\cc\ |\ l(\phi^-X)>0\}$). We call $l$ an {\it additive function} if it is left-right additive. Put $l(a):=l(X)-l(Y)$ for $a\in\cc(X,Y)$.

\vskip1em{\bf\XACB\ Examples }(1) Let $\Lambda$ be an artin algebra and $\cc:=\mod\Lambda$. Then $\cc$ forms a strict $\tau$-category with $\pp\cc=\ind(\pr\Lambda)$ and $\ii\cc=\ind(\inj\Lambda)$. Moreover, $(X]$ gives an almost split sequence for any $X\in\ind\cc-\pp\cc$ [ARS], and $l(X):=\length_\Lambda X$ gives an additive function with $l^\pm=\ind^\pm_1\cc$. More generally, since any contravariantly finite torsionfree class $\cc$ over $\Lambda$ has almost split sequences [AS], $\cc$ forms a strict $\tau$-category.

(2) Let $\Gamma$ be a semiperfect ring and $\cc:=\pr\Gamma$. Then $\cc$ forms a strict $\tau$-category if and only if $\gl\Gamma\le2$ and any simple $\Gamma$ or $\Gamma^{op}$-module $L$ with $\pd L=2$ satisfies that $\grade{L}=2$ and $\ext_\Gamma^2(L,\Gamma)$ is a simple $\Gamma^{op}$ or $\Gamma$-module.
In this case, $\pp\cc=\{ P\in\ind\cc\ |\ \pd\top P\le1\}$ and $\ii\cc=\{ P\in\ind\cc\ |\ \pd\top\widehat{P}\le1\}$ hold.

(3) Let $\qq$ be a $\tau$-species (=modulated translation quiver in [IT2]) and $\cc$ the mesh category of $\qq$ [I3;8.3,8.4]. Then $\cc$ is a (not necessarily strict) $\tau$-category. Thus we obtain a bijection between isomorphism classes of $\tau$-species and equivalence classes of completely graded $\tau$-categories [I3;10.1] by taking mesh categories. This structure theorem of completely graded $\tau$-categories was a strong motivation of the introduction of $\tau$-categories in [I3].

\vskip1em{\bf\XACC\ }[I3;4.1,7.2] Let $\cc$ be a $\tau$-category.
For $X=\sum_{Y\in Q}a_YY\in\zzz\ind\cc$, put $X_+:=\sum_{Y\in Q,a_Y>0}a_YY\in\nnn\ind\cc$. Define a map $\theta^+_n:\nnn\ind\cc\rightarrow\nnn\ind\cc$ ($n\ge0$) by $\theta^+_0:=1_{\nnn Q}$, $\theta^+_1:=\theta^+$ and $\theta^+_nX:=(\theta^+\theta^+_{n-1}X-\tau^+\theta^+_{n-2}X)_+$ for $n\ge2$. Then $\theta^+_n$ becomes a monoid momorphism, and for any $X\in\cc$, there exists the following commutative diagram such that $\cc(\ ,\tau^+\theta^+_{n-1}X)\stackrel{\cdot a_n}{\longrightarrow}\cc(\ ,\theta^+_nX)\rightarrow\jac_{\cc}^n(\ ,X)\rightarrow0$ is exact and $a_n$ is in $\jac_{\cc}$ for any $n\ge0$.
\[\begin{diag}
\cdots&\RA{}&\tau^+\theta^+_3X&\RA{}&\tau^+\theta^+_2X&\RA{}&\tau^+\theta^+_1X&\RA{}&\tau^+X&\RA{}&0\\
&&\downarrow^{a_4}&&\downarrow^{a_3}&&\downarrow^{a_2}&&\downarrow^{a_1}&&\downarrow^{a_0}\\
\cdots&\RA{}&\theta^+_4X&\RA{}&\theta^+_3X&\RA{}&\theta^+_2X&\RA{}&\theta^+_1X&\RA{}&X\end{diag}\]

\vskip0em{\bf\XAD\ }
Let $\cc$ be a $\tau$-category and $A,B\in\ind\cc$. We say that {\it $(A,B)$ is a Nakayama pair} if there exists the following commutative diagram for some $n\ge0$ such that $(X_{i}\stackrel{(a_{i}\ g_{i})}{\longrightarrow}Y_{i}\oplus X_{i-1}\stackrel{{-f_{i}\choose a_{i-1}}}{\longrightarrow}Y_{i-1})$ is a $\tau$-sequence for any $i$ ($0<i\le n$).
\[\begin{diag}
&X_n&&\RA{g_n}&X_{n-1}&\RA{g_{n-1}}&&\RA{g_2}&X_1&\RA{g_1}&X_0&=&A\\
&\downarrow&{}^{a_n=\mu^+_B}&&\downarrow^{a_{n-1}}&&\cdots&&\downarrow^{a_1}&&&\downarrow&{}^{a_0=\mu^-_A}\\
B&=&Y_n&\RA{f_n}&Y_{n-1}&\RA{f_{n-1}}&&\RA{f_2}&Y_1&\RA{f_1}&&Y_0&\end{diag}\]

Define $\eta^+_i\in\endm_{\zzz}(\zzz\ind\cc)$ ($i\ge0$) by $\eta^+_0:=\theta^-$, $\eta^+_1:=\theta^+\circ\theta^--1_{\zzz\ind\cc}$ and $\eta^+_i:=\theta^+\circ\eta^+_{i-1}-\tau^+\circ\eta^+_{i-2}$ for $i\ge2$.
Then $Y_i=\eta^+_iA$ and $X_{i+1}=\tau^+\eta^+_iA$ holds immediately.
In particular, $B$ is uniquely determined by $A$, and vice versa.
We write $B=\nak^-(A)$ and $A=\nak^+(B)$. Note that any right (or left) additive function $l$ satisfies $l(a_0)=l(a_1)=\cdots=l(a_n)$.

\vskip1em{\bf\XADA\ Example }
(1) Let $\Lambda$ be an order over a complete discrete valuation ring $R$ and $\cc$ the category of $\Lambda$-lattices [CR]. Then $\cc$ forms a $\tau$-category. If $\Lambda$ is representation-finite, then $(A,B)$ is a Nakayama pair for any $B\in\ind(\pr\Lambda)$ and $A:=\hom_R(\widehat{B},R)$ by [I4]3.3.

(2) Let $\Lambda$ be a representation-finite artin algebra and $\cc:=\mod\Lambda$. Let $B\in\ind(\pr\Lambda)$, $A:=(\widehat{B})^*\in\ind(\inj\Lambda)$ and $X:=\soc A=\top B\in\ind(\ssim\Lambda)$. In \XCCA, we will show that $(A,\tau^-X)$ is a Nakayama pair if $A$ is not simple, and $(\tau^+X,B)$ is a Nakayama pair if $B$ is not simple.

\vskip1em{\bf\XADB\ }We collect basic results on Nakayama pairs, where we refer [I4;8.1] to (1), [I3;6.4] to (2) and [I4;2.2] to (3).
Let $\cc$ be a $\tau$-category.

(1) $(A,B)$ is a Nakayama pair if and only if there exists $n\ge0$ such that $\eta^+_iA\in\nnn(\ind\cc-\pp\cc)$ for any $i$ ($0\le i<n$), $\eta^+_nA=B$ and $\eta^+_{n+1}A=0$.

(2) Assume $\bigcap_{n\ge0}\jac_{\cc}^n=0$. If $\mu^-_A$ is not a monomorphism for $A\in\ind\cc-\ind^-_0\cc$, then $(A,B)$ is a Nakayama pair for some $B\in\ind\cc-\pp\cc$.

(3) Assume that $\cc=\pr\Gamma$ for a semiperfect ring $\Gamma$. For $A,B\in\ind\cc$, let $L:=\top A$ be a simple $\Gamma$-module and $M:=\top\widehat{B}$ a simple $\Gamma^{op}$-module. Then $(A,B)$ is a Nakayama pair if and only if $\tr L$ has finite length with the socle $M$ and $\sgrade{(\tr L)/M}\ge2$ if and only if $\tr M$ has finite length with the socle $L$ and $\sgrade{(\tr M)/L}\ge2$.

\vskip1em{\bf\XAE\ } $\tau$-categories were defined `locally' by the properties of simple modules over the category [I3]. On the other hand, Rump [R1,2] introduced the concept of almost abelian categories, which is given `globally' by the properties of kernels and cokernels.

An additive category is called {\it preabelian} if any morphism has a kernel and a cokernel. A preabelian category is called {\it almost abelian} if kernels are stable under pushout and cokernels are stable under pullback. An almost abelian category is called {\it integral} if monomorphisms are stable under pushout and epimorphisms are stable under pullback. Rump has shown that they are closely related to tilting theory.

\vskip1em{\bf\XB\ Repesentation theoretic realization of artinian strict $\tau$-categories}

\vskip1em
{\bf\XBA\ Theorem }{\it Let $\Gamma$ be an artin algebra and $\cc:=\pr\Gamma$.
Then the following conditions are equivalent.

(1) $\Gamma$ satisfies $\gl\Gamma\le2$ and the $(2,2)$ and $(2,2)^{op}$-condition (\S\XAA).

(2) There exists an artin algebra $\Lambda$ and a (faithful) torsion theory $(\tt,\ff)$ on $\mod\Lambda$ such that $\cc$ is equivalent to $\ff$ (\S\XAB).

(3) $\cc$ is a strict $\tau$-category (\S\XAC).

(4) $\cc$ is a $\tau$-category with a right additive function (\S\XACA).

(5) $\cc$ is an almost abelian category (\S\XAE).

($i$)$^{op}$ Opposite side version of ($i$) ($1\le i\le 5$).}

\vskip1em{\bf\XBAA\ Remark }
(1) The equivalence of \XBA(2) and (2)$^{op}$ follows from the classical cotilting theorem (see \XABA(2)), the equivalence of \XBA(2) and (3) is a special case of [R3;Theorem 1], and (2)$\Rightarrow$(3) follows from the remark in \XACB(1).

(2) There exists an artin algebra $\Gamma$ with $\gl\Gamma\le2$ such that $\Gamma$ satisfies exactly one of the $(2,2)$ and $(2,2)^{op}$-condition. For example, such an algebra $\Gamma$ is given by the quiver {\scriptsize$\bullet\stackrel{a}{\rightarrow}\bullet\stackrel{b}{\rightarrow}\bullet\leftarrow\bullet$} with the relation $ba=0$.

\vskip1em{\bf\XBB\ Definition }Let $\Gamma$ be an artin algebra and $I^i:=I^i_\Gamma(\Gamma)$.
We call a functor $\ppp:\pr\Gamma\rightarrow\mod\Lambda$ a ({\it representation theoretic}) {\it realization} of $\Gamma$ if $\Lambda$ is an artin algebra, $\ppp$ is full faithful and $\Lambda\in\ppp(\pr\Gamma)$. This is equivalent to that, there exists $Q\in\pr\Gamma$ such that $\Lambda=\endm_\Gamma(Q)$, $\ppp$ is isomorphic to $\hom_\Gamma(Q,\ )$ and $I^0\oplus I^1\in\add(\widehat{Q})^*$. We call a realization $\ppp=\hom_\Gamma(Q,\ )$ {\it minimal} if $\add(I^0\oplus I^1)=\add(\widehat{Q})^*$ holds. A minimal realization of $\Gamma$ is unique up to Morita-equivalence. We sometimes regard $\ppp=\hom_\Gamma(Q,\ )$ as a functor $\mod\Gamma\rightarrow\mod\Lambda$.


\vskip1em{\bf\XBBA\ }
Let $\Gamma$ be a noetherian ring and $I^i:=I^i_\Gamma(\Gamma)$. For $Q\in\pr\Gamma$, put $\Lambda:=\endm_\Gamma(Q)$ and $\ppp:=\hom_\Gamma(Q,\ ):\mod\Gamma\rightarrow\mod\Lambda$. Then $\ppp$ induces an equivalence $\add Q\rightarrow\pr\Lambda$, and the conditions (i)--(iv) below are equivalent.

(i) $\grade{X}\ge2$ (resp. $\grade{X}\ge1$) holds for any $X\in\mod\Gamma$ with $\ppp X=0$.

(ii) $\ppp$ is full faithful (resp. faithful) on $\pr\Gamma$.

(iii) $\ppp_{X,Y}$ is bijective (resp. injective) for any $X\in\mod\Gamma$ and $Y\in\pr\Gamma$.

(iv) $I^0\oplus I^1\in\add(\widehat{Q})^*$ (resp. $I^0\in\add(\widehat{Q})^*$).

\vskip1em{\sc Proof }
We only prove the assertion for `full faithful'. (iii)$\Rightarrow$(ii) and (i)$\Leftrightarrow$(iv) are clear (see \XCB(1)).

(ii)$\Rightarrow$(i) Let $P_2\rightarrow P_1\rightarrow P_0\rightarrow X\rightarrow0$ be a projective resolution. Then $\ppp P_2\rightarrow\ppp P_1\rightarrow\ppp P_0\rightarrow0$ is exact. Since $\ppp$ is full faithful on $\pr\Gamma$, we obtain an exact sequence $\widehat{P}_2\leftarrow\widehat{P}_1\leftarrow\widehat{P}_0\leftarrow0$ by taking $\hom_{\Lambda}(\ ,\ppp\Gamma)$. Thus $\grade{X}\ge2$ holds.

(i)$\Rightarrow$(iii) We can take a complex ${\bf A}:Q^m\stackrel{f_1}{\rightarrow}Q^n\stackrel{f_0}{\rightarrow}X\rightarrow0$ such that $\Lambda^m\stackrel{\ppp f_1}{\longrightarrow}\Lambda^n\stackrel{\ppp f_0}{\longrightarrow}\ppp X\rightarrow0$ is exact. We obtain an exact sequence $\ppp \Gamma^m\stackrel{\ppp f_1\cdot}{\longleftarrow}\ppp \Gamma^n\stackrel{\ppp f_0\cdot}{\longleftarrow}\hom_{\Lambda}(\ppp X,\ppp\Gamma)\leftarrow0$ by taking $\hom_{\Lambda}(\ ,\ppp\Gamma)$. On the other hand, since $\ppp(\Cok f_j)=0$ ($j=0,1$) holds, we obtain $\grade{\Cok f_j}\ge2$. Taking $\hom_{\Gamma}(\ ,\Gamma)$ for ${\bf A}$, we obtain an exact sequence $\ppp\Gamma^m\stackrel{\ppp f_1\cdot}{\longleftarrow}\ppp\Gamma^n\stackrel{\ppp f_0\cdot}{\longleftarrow}\hom_{\Gamma}(X,\Gamma)\leftarrow0$.\rule{5pt}{10pt}

\vskip1em{\bf\XBBB\ Proof of \XBB\ }
Assume that $Q\in\pr\Gamma$ satisfies $I^0\oplus I^1\in\add(\widehat{Q})^*$ and put $\Lambda:=\endm_\Gamma(Q)$. Then $\hom_\Gamma(Q,\ )$ is a realization of $\Gamma$ by \XBBA. Conversely, let $\ppp$ be a realization of $\Gamma$. Take $Q\in\pr\Gamma$ such that $\ppp Q=\Lambda$. Since $\hom_\Gamma(Q,\ )\stackrel{\ppp}{\stackrel{\sim}{\rightarrow}}\hom_\Lambda(\ppp Q,\ppp(\ ))=\ppp$ holds, $\ppp$ is isomorphic to $\hom_\Gamma(Q,\ )$. Moreover, $I^0\oplus I^1\in\add(\widehat{Q})^*$ holds by \XBBA.\rule{5pt}{10pt}

\vskip1em{\bf\XBBC\ }
(1) Let $\Lambda$ be an artin algebra, $\cc$ a full subcategory of $\mod\Lambda$ with an additive generator $M$ and $\Gamma:=\endm_\Lambda(M)$. Then the functors $\qqq:=\hom_\Lambda(M,\ ):\mod\Lambda\rightarrow\mod\Gamma$ and $\rrr:=\hom_\Lambda(\ ,M):\mod\Lambda\rightarrow\mod\Gamma^{op}$ induce equivalences $\cc\rightarrow\pr\Gamma$ and $\cc\rightarrow\pr\Gamma^{op}$ such that $\hom_\Gamma(\ ,\Gamma)\circ\qqq=\rrr$ and $\hom_\Gamma(\ ,\Gamma)\circ\rrr=\qqq$ hold on $\cc$.

(2) Let $\ppp:\pr\Gamma\rightarrow\mod\Lambda$ be a realization of $\Gamma$. Put $\cc:=\ppp(\pr\Gamma)$ and $M:=\ppp\Gamma$. Then $\qqq:\cc\rightarrow\pr\Gamma$ in (1) gives a quasi-inverse of $\ppp:\pr\Gamma\rightarrow\cc$.

\vskip1em{\bf\XBC\ Lemma }{\it
Let $\Gamma$ be an artin algebra and $\ppp:\pr\Gamma\rightarrow\mod\Lambda$ a realization. Assume that $\Gamma$ satisfies $\gl\Gamma\le2$ and the $(2,2)^{op}$-condition. Then $\ppp(\pr\Gamma)$ is closed under kernels and extensions in $\mod\Lambda$.}

\vskip1em{\sc Proof }
Take $Q\in\pr\Gamma$ in \XBB\ and extend $\ppp=\hom_\Gamma(Q,\ )$ to $\mod\Gamma\rightarrow\mod\Lambda$. Using $\gl\Gamma\le2$, we can easily show that $\pr\Gamma$ is closed under kernels.

(i) Let ${\bf A}:0\stackrel{}{\rightarrow}\ppp P^\prime\stackrel{g}{\rightarrow}X\stackrel{f}{\rightarrow}\ppp P\rightarrow0$ be an exact sequence in $\mod\Lambda$ with $P,P^\prime\in\pr\Gamma$.
We will show that there exists an exact sequence ${\bf B}:0\stackrel{}{\rightarrow}P^\prime\stackrel{}{\rightarrow}L\stackrel{}{\rightarrow}P\rightarrow M\rightarrow0$ in $\mod\Gamma$
such that ${\bf A}$ is isomorphic to $\ppp{\bf B}$ as a complex.

Take a surjection $d\in\hom_\Lambda(\ppp P_1,X)$ with $P_1\in\add Q$ by \XBBA, and take $a\in\hom_{\Gamma}(P_1,P)$ such that $df=\ppp a$. Taking an exact sequence ${\bf C}:0\rightarrow\Omega^2M\stackrel{b}{\rightarrow}P_1\stackrel{a}{\rightarrow}P\rightarrow M\rightarrow0$, we obtain the following commutative diagram.
\[\begin{diag}
{\bf A}:&0&\RA{}&\ppp P^\prime&\RA{g}&X&\RA{f}&\ppp P&\RA{}&0\\
&&&\uparrow^{e}&&\uparrow^{d}&&\parallel\\
\ppp{\bf C}:&0&\RA{}&\ppp\Omega^2M&\RA{\ppp b}&\ppp P_1&\RA{\ppp a}&\ppp P&\RA{}&0\end{diag}\]

Take $c\in\hom_{\Gamma}(\Omega^2M,P^\prime)$ such that $\ppp c=e$ by \XBBA, and define ${\bf B}$ by the following push-out diagram.
\[\begin{diag}
{\bf B}:&0&\RA{}&P^\prime&\RA{}&L&\RA{}&P&\RA{}&M&\RA{}&0\\
&&&\uparrow^{c}&&\uparrow^{}&&\parallel&&\parallel\\
{\bf C}:&0&\RA{}&\Omega^2M&\RA{b}&P_1&\RA{a}&P&\RA{}&M&\RA{}&0
\end{diag}\]

Since both of ${\bf A}$ and $\ppp{\bf B}$ are given by the push-out of $\ppp{\bf C}$ by $c=\ppp h$, the complexes ${\bf A}$ and $\ppp{\bf B}$ are isomorphic.

(ii) To show the lemma, take the complex ${\bf B}$ in (i).
Since $\grade{M}\ge2$ holds by $\ppp M=0$, we have an exact sequence $\ext^2_\Gamma(M,\Gamma)\leftarrow\widehat{P}^\prime\leftarrow\widehat{L}\leftarrow\widehat{P}\leftarrow0$ by taking $\widehat{(\ )}$. Since $\Gamma$ satisfies the $(2,2)^{op}$-condition, we obtain the upper exact sequence of the following commutative diagram by taking $\widehat{(\ )}$ again.
\[\begin{diag}
&0&\RA{}&\widehat{\widehat{P}}{}^\prime&\RA{}&\widehat{\widehat{L}}&\RA{}&\widehat{\widehat{P}}&\\
&&&\parallel&&\uparrow&&\parallel\\
{\bf B}:&0&\RA{}&P^\prime&\RA{}&L&\RA{}&P&\RA{}&M&\RA{}&0
\end{diag}\]

Taking the mapping cone, we obtain an exact sequence $0\rightarrow L\rightarrow\widehat{\widehat{L}}\rightarrow M$. Thus $X=\ppp L=\ppp\widehat{\widehat{L}}$ holds. Since $\gl\Gamma\le2$ holds, we obtain $\widehat{L}\in\pr\Gamma^{op}$. Thus $X\in\ppp(\pr\Gamma)$.\rule{5pt}{10pt}

\vskip1em{\bf\XBD\ Proof of \XBA\ }
(1)$\Leftrightarrow$(1)$^{op}$ is clear, and (1)$\Leftrightarrow$(3) holds by [I5]6.3.

(1)$\Rightarrow$(2) Let $\ppp=\hom_\Gamma(Q,\ )$ be a minimal realization of $\Gamma$ which we extend to $\mod\Gamma\rightarrow\mod\Lambda$, and $\ff:=\ppp(\pr\Gamma)$ a full subcategory of $\mod\Lambda$. Then $\ff$ is closed under extensions by \XBC. We will show that $\ff$ is closed under submodules. Fix any $P\in\pr\Gamma$ and an injection $f\in\hom_\Lambda(X,\ppp P)$. Take a surjection $g\in\hom_\Lambda(\ppp P_1,X)$ with $P_1\in\add Q$ by \XBBA, and take $a\in\hom_\Gamma(P_1,P)$ such that $fg=\ppp a$. Then $X=\ppp L$ holds for $L:=\Im a$. The set $\{ M\in\mod\Gamma\ |\ L\subseteq M\subseteq P,\ \ppp(M/L)=0\}$ has a unique maximal element, which we denote by $M$. Then $\ppp M=X$ and $\soc(P/M)\in\add\top Q=\add\soc(I^0_\Gamma(\Gamma)\oplus I^1_\Gamma(\Gamma))$ hold. Since $\Gamma$ satisfies the $(2,2)$-condition, an injective hull $I$ of $\soc(P/M)$ satisfies $\pd I\le 1$. Since $I$ gives an injective hull of $P/M$, we obtain $\pd P/M\le1$ by $\gl\Gamma\le2$. Thus $X=\ppp M$ and $M\in\pr\Gamma$ hold, and $\ff$ is a faithful torsionfree class.


(2)$\Rightarrow$(4) We will show that $l(X):=\length_\Lambda X$ gives a right additive function. For any $X\in\cc$, $(X]$ gives an exact sequence $0\rightarrow\tau^+X\stackrel{\nu^+_X}{\rightarrow}\theta^+X\stackrel{\mu^+_X}{\rightarrow}X$ in $\mod\Lambda$ since the kernel of $\mu^+_X$ in $\mod\Lambda$ is contained in $\ff$ and thus coincides with $\tau^+X$. Hence $l(\phi^+X)=l(X)-l(\theta^+X)+l(\tau^+X)\ge0$ holds. We only have to show that $\mu^+_X$ is surjective for any $X\in\ind\cc-\pp\cc$. Otherwise, the inclusion $Y:=\Im\mu^+_X\stackrel{a}{\rightarrow}X$ induces an isomorphism $\cc(\ ,Y)\stackrel{\cdot a}{\rightarrow}\jac_{\cc}(\ ,X)$ with $Y\in\ff$. Thus $0\rightarrow Y\stackrel{a}{\rightarrow}X$ gives $(X]$, a contradiction to $X\notin\pp\cc$.

(4)$\Rightarrow$(3) Let $l$ be a right additive function.
For any $X\in\ind\cc-\pp\cc$, we only have to show that $\nu^+_X$ is a monomorphism. Otherwise, there exists $Y\in\ind\cc-\pp\cc$ such that $(\tau^+X,Y)$ is a Nakayama pair by \XADB(2).
Then $0>-l(X)=l(\nu^+_X)=l(\mu^+_Y)=l(\tau^+Y)>0$ holds by \XAD, a contradiction.

(2)$\Leftrightarrow$(5) See [R2;Theorem 1].\rule{5pt}{10pt}

\vskip1em{\bf\XBE\ }We call an artin algebra $\Gamma$ a {\it strict $\tau$-algebra} if it satisfies the equivalent conditions in \XBA.
We denote by $\gd\Gamma:=\bigoplus_{n\ge0}J_\Gamma^n/J_\Gamma^{n+1}$ the associated graded algebra.
Radical Layers Theorem of Igusa-Todorov ([IT1][BG]), which is one of the most important theorems in the representation theory of algebras, was proved for arbitrary artin algebras and even for $\tau$-categories in [I3;4.2]. Consequently we obtain \XBEA\ below, which implies the following theorem immediately [I3;5.2].

\vskip1em{\bf Theorem }{\it
Let $\Gamma$ be a strict $\tau$-algebra. Then $\Gamma$ is strongly quasi-Koszul in the sense of Green-Martinez [GM;\S5], and $\gd\Gamma$ is a strict $\tau$-algebra again.}

\vskip1em{\bf\XBEA\ Lemma }{\it Let $\Gamma$ be a strict $\tau$-algebra and $0\rightarrow P_2\stackrel{g}{\rightarrow}P_1\stackrel{f}{\rightarrow}P_0\rightarrow L\rightarrow0$ a minimal projective resolution of a simple $\Gamma$-module $L$. Then $0\rightarrow J_\Gamma^{i-1}P_2\stackrel{g}{\rightarrow}J_\Gamma^{i}P_1\stackrel{f}{\rightarrow}J_\Gamma^{i+1}P_0\rightarrow0$ is exact for any $i\ge0$, where we put $J_\Gamma^{-1}:=\Gamma$}

\vskip1em{\bf\XC\ Auslander-regular artin algebra with global dimension two}

In this section, we study several variations of our theorem \XBA.

\vskip1em{\bf\XCA\ Theorem }{\it Let $\Gamma$ be an artin algebra and $\cc:=\pr\Gamma$. Then the following conditions are equivalent.

(1) $\Gamma$ is an Auslander-regular ring with $\gl\Gamma\le2$ (\S\XAA).

(2) There exists an artin algebra $\Lambda$ and a (faithful) hereditary torsion theory $(\tt,\ff)$ on $\mod\Lambda$ such that $\cc$ is equivalent to $\ff$ (\S\XAB).

(3) $\cc$ is a strict $\tau$-category and $\nak^-$ gives a map $\ii\cc-\ind^-_0\cc\rightarrow\ind\cc$ (\S\XAD).

(4) $\cc$ is a $\tau$-category with an additive function (\S\XACA).

(5) $\cc$ is an integral almost abelian category (\S\XAE).

($i$)$^{op}$ Opposite side version of ($i$) ($1\le i\le 5$).}

\vskip1em{\bf\XCAA\ Lemma }{\it
Let $\Lambda$ and $\Gamma$ be artin algebras and $(\tt,\ff)$ a torsion theory
on $\mod\Lambda$. Assume that $\pr\Gamma$ is equivalent to $\ff$.
Then (1) and (2) below are equivalent, and (3) implies them.
If $(\tt,\ff)$ is faithful, then (1)--(3) are equivalent.

(1) $\Gamma$ is an Auslander-regular ring with $\gl\Gamma\le2$.

(2) If $f\in\hom_\Lambda(Y,X)$ is a surjection with $Y\in\ff$ and $X\in\tt$, then $f(\soc Y)=0$.

(3) $(\tt,\ff)$ is hereditary.}

\vskip1em{\sc Proof }
$\Gamma$ satisfies $\gl\Gamma\le2$ and the $(2,2)$ and $(2,2)^{op}$-condition, and $\ff$ forms a $\tau$-category by \XBA. We use the notations in \XBBC\ and \XAB(iv).

(2)$\Rightarrow$(1) To show that $\Gamma$ satisfies the $(1,1)^{op}$-condition, we will show $\grade{M}>0$ for any simple $\Gamma$-module $L$ and a submodule $M$ of $\ext^1_\Gamma(L,\Gamma)$.
Since $\pd L=2$ implies $\grade{L}=2$ by \XACB(2), we can assume $\pd L=1$.
Take projective resolutions $0\rightarrow\qqq Y\stackrel{\qqq f}{\longrightarrow}\qqq X\rightarrow L\rightarrow0$ and $0\leftarrow M\stackrel{a}{\leftarrow}\rrr W$. Then $0\leftarrow\ext^1_\Gamma(L,\Gamma)\leftarrow\rrr Y\stackrel{\rrr f}{\longleftarrow}\rrr X$ is exact, and $\phi$ lifts to $\rrr d\in\hom_\Gamma(\rrr W,\rrr Y)$. Take an exact sequence $Y\stackrel{(f\ d)}{\longrightarrow}X\oplus W\stackrel{{g\choose e}}{\rightarrow}V\rightarrow0$.
Then $0\leftarrow M\stackrel{a}{\leftarrow}\rrr W\stackrel{\rrr(e\alpha_V)}{\longleftarrow}\rrr\circ\fff V$ gives a projective resolution for $e\alpha_V\in\hom_\Lambda(W,\fff V)$. Since $0\rightarrow\widehat{M}\rightarrow\qqq W\stackrel{\qqq(e\alpha_V)}{\longrightarrow}\qqq\circ\fff V$ is exact, we only have to show that $e\alpha_V$ is injective. Put $U:={g\choose e}^{-1}(\ttt V)\subseteq X\oplus W$. Then ${g\choose e}(\soc U)=0$ holds by (2). Since $f$ and $e$ are injective, $W\cap\soc U=0$ holds. Thus $W\cap U=0$ holds, and we obtain the assertion.

(1)$\Rightarrow$(2) Let $0\rightarrow Z\stackrel{g}{\rightarrow}Y
\stackrel{f}{\rightarrow}X\rightarrow0$ be an exact sequence such that $Y\in\ff$, $X\in\tt$ and $\soc Y\ {\not\subseteq}\ Z$.
Then there exists an injection ${g\choose h}\in\hom_\Lambda(Z\oplus W,Y)$ with $W\neq0$.
Define $L\in\mod\Gamma$ by an exact sequence $0\rightarrow\qqq(Z\oplus W)\stackrel{\qqq{g\choose h}}{\longrightarrow}\qqq Y\rightarrow L\rightarrow0$.
Then $0\leftarrow\ext^1_\Gamma(L,\Gamma)\stackrel{a}{\leftarrow}\rrr(Z\oplus W)\stackrel{\rrr{g\choose h}}{\longleftarrow}\rrr Y$ is exact.
Since $\rrr g$ is injective, $a$ restricts to an injection $\rrr W\rightarrow\ext^1_\Gamma(L,\Gamma)$. Thus $\qqq W=\hom_\Gamma(\rrr W,\Gamma)=0$, a contradiction.

(3)$\Rightarrow$(2) For any simple submodule $Z$ of $Y$, $Z\in\ff$ and $f(Z)\in\tt$ imply $f(Z)=0$.

(2)+($\Lambda\in\ff$)$\Rightarrow$(3) Let $X$ be a submodule of $Y\in\tt$ and $f\in\hom_\Lambda(\fff X,Y/\ttt X)$ a natural injection. Take a surjection ${g\choose f}\in\hom_\Lambda(P\oplus\fff X,Y/\ttt X)$ with $P\in\pr\Lambda$. By $P\oplus\fff X\in\ff$ and $Y/\ttt X\in\tt$, we obtain $f(\soc\fff X)=0$. Thus $\soc\fff X=0$ and $X\in\tt$.\rule{5pt}{10pt}

\vskip1em{\bf\XCAB\ Proof of \XCA\ }(4)$\Leftrightarrow$(4)$^{op}$ holds clearly.


(3)$\Rightarrow$(1) We use the notations in \XBBC. To show that $\Gamma$ satisfies the $(1,1)^{op}$-condition, we will show $\sgrade{\ext^1_\Gamma(L,\Gamma)}\ge1$ for any simple $\Gamma$-module $L$ with $\pd L=1$.
Put $L=\top\rrr A$ for $A\in\ii\cc-\ind^-_0\cc$, $B:=\nak^-(A)\in\ind\cc$ and $M:=\top\qqq B$.
Then we have an exact sequence $0\rightarrow M\rightarrow\ext^1_\Gamma(L,\Gamma)\rightarrow N\rightarrow0$, where $\sgrade{N}\ge2$ holds by \XADB(3).
Since $\grade{\ext^1_\Gamma(L,\Gamma)}\ge1$ holds by $\pd L=1$, we obtain $\sgrade{M}=\grade{M}\ge1$ by taking $\hom_\Gamma(\ ,\Gamma)$. Thus $\sgrade{\ext^1_\Gamma(L,\Gamma)}\ge1$ holds.

(1)$\Rightarrow$(2) Immediate from \XBA\ and \XCAA.

(2)$\Rightarrow$(4) By \XAB, we can put $\ff={\bf S}^\perp$ and $\tt={}^\perp\ff$ for ${\bf S}\subseteq\ind(\ssim\Lambda)$. For $X\in\mod\Lambda$, we denote by $l(X)$ the number of its composition factors which is {\it not} in ${\bf S}$. We will show that $l$ gives an additive function. By the argument in the proof of \XBA(2)$\Rightarrow$(4), $l$ is right additive. For any $X\in\ii\cc$, take an exact sequence $X\stackrel{\mu^-_X}{\rightarrow}\theta^-X\rightarrow Y\rightarrow0$ in $\mod\Lambda$. Since $\mu^-_X$ is an epimorhism in $\cc$, we obtain $Y\in\tt$ and $l(Y)=0$. Thus $l(X)\ge l(\theta^-X)$ holds for any $X\in\ii\cc$, and $l$ is additive.

(4)$\Rightarrow$(3) $\cc$ is strict by \XBA. Let $l$ be an additive function. In the proof below, we have to use concepts in [I3]. Take $X\in\ii\cc-\ind^-_0\cc$ and let ${\bf a}=(a_i)_{0\le i}$ be the right ladder of $\mu^-_X$ for $a_i\in\cc(X_i,Y_i)$.
If ${\bf a}$ is not essential, then $(X,Y)$ is a Nakayama pair for some $Y\in\ind\cc-\pp\cc$ by [I3;6.4]. Thus we assume that ${\bf a}$ is essential.
Take a maximal number $n$ such that $Y_n\neq0$. Then $Y_n\in\add(\pp\cc)$ and $a_n=\mu^+_{Y_n}$ hold. Let ${\bf c}=(c_i)_{0\le i\le n}$ be the left ladder of $a_n$ for $c_i\in\cc(A_i,B_i)$. Since $A_i$ has no direct summands in $\ii\cc$ for any $i$ ($0\le i<n$) by [I3;6.3.1(1)], ${\bf c}$ is invertible by [I3;6.2.1]. Since $l(c_i)=l(a_n)\le0$ holds, $B_i=0$ implies $A_i=0$. Hence $A_i\neq0$ holds for any $i$ ($0\le i\le n$) inductively. Since $(a_i)_{0\le i\le n}$ is invertible by [I3;6.3.1(2)(i)], $(X,Y_n)$ is a Nakayama pair and $\nak^-(X)=Y_n$.

(2)$\Rightarrow$(5) follows from [R2;Lemma 6], and (5)$\Rightarrow$(1) follows from a quite similar argument as in the proof of \XCAA(2)$\Rightarrow$(1).\rule{5pt}{10pt}


\vskip1em{\bf\XCB\ }
Let $\Gamma$ be an artin algebra with $I^i:=I^i_\Gamma(\Gamma)$.

(1) The bijection $\soc:\ind(\inj\Gamma)\rightarrow\ind(\ssim\Gamma)$ induces the maps below. The first and second maps are bijective, and so is the third map if $\gl\Gamma<\infty$.
{\scriptsize\begin{eqnarray*}
\ind(\add I^n)&\rightarrow&\{ L\in\ind(\ssim\Gamma)\ |\ \ext^n_\Gamma(L,\Gamma)\neq0\}\\
\ind(\add I^n)-\bigcup_{i<n}\ind(\add I^i)&\rightarrow&\{ L\in\ind(\ssim\Gamma)\ |\ \grade{L}=n\}\\
\ind(\add I^n)-\bigcup_{i>n}\ind(\add I^i)&\rightarrow&\{ L\in\ind(\ssim\Gamma)\ |\ \pd L=n\}
\end{eqnarray*}}

\vskip-1em
(2) Assume that the conditions in \XBA\ are satisfied and $(\tt,\ff)$ is faithful. Let $\ppp:\pr\Gamma=\cc\stackrel{\sim}{\rightarrow}\ff\subset\mod\Lambda$ be the composition. Then $\ind(\add(I^0\oplus I^1))\cap\ind(\add I^2)=\emptyset$ holds, and $\ppp$ is a minimal realization of $\Gamma$ given by $\ppp=\hom_\Gamma(Q,\ )$ for some $Q\in\pr\Gamma$. We have the bijections below, where the middle map is given by the projective cover.
{\scriptsize\[\begin{diag}
\ind(\inj\Gamma)&\stackrel{\soc}{\longrightarrow}&\ind(\ssim\Gamma)&\stackrel{}{\longrightarrow}&\ind(\pr\Gamma)=\ind\cc&\stackrel{\ppp}{\longrightarrow}&\ind(\ff)\\
\cup&&\cup&&\cup&&\cup\\
\ind(\add(I^0\oplus I^1))&\longrightarrow&\{ L\in\ind(\ssim\Gamma)\ |\ \pd L\le1\}&\longrightarrow&\ind(\add Q)=\pp\cc&\longrightarrow&\ind(\pr\Lambda)\\
\cup&&\cup&&\cup&&\cup\\
\ind(\add(I^0\oplus I^1))-\ind(\add I^1)&\longrightarrow&\{ L\in\ind(\ssim\Gamma)\ |\ \pd L=0\}&\longrightarrow&\ind^+_0\cc&\longrightarrow&\ind(\pr\Lambda)\cap\ind(\ssim\Lambda)
\end{diag}\]}

\vskip-1em
(3) In (2), assume that the conditions in \XCA\ are satisfied and $l$ is an additive function of $\cc$. Then the bijections in (2) induce the bijections (i), and the equalities (ii) hold:
{\scriptsize\begin{eqnarray*}
&\ind(\add I^0)\longrightarrow\{ L\in\ind(\ssim\Gamma)\ |\ \grade{L}=0\}\longrightarrow l^+\longrightarrow\{ X\in\ind(\pr\Lambda)\ |\ \top X\in\ff\}&\mbox{(i)}\\
&l^+=\ind^+_0\cc\cup\{ X\in\pp\cc-\ind^+_0\cc\ |\ \nak^+(X)\in\ind\cc-\ii\cc\}=\{ X\in\pp\cc\ |\ \mu^+_X\mbox{ is not an epimorphism}\}&\mbox{(ii)}
\end{eqnarray*}}

\vskip-1em{\sc Proof }
(1) Taking $\hom_\Gamma(L,\ )$ for a minimal injective resolution $0\rightarrow \Gamma\stackrel{}{\rightarrow}I^0\stackrel{}{\rightarrow}I^1\stackrel{}{\rightarrow}\cdots$, we obtain the first bijection, which implies others.

(2) The former assertion follows from \XACB(2), and $\soc$ induces the left bijections by (1). Since $\ppp$ is a realization by $\Lambda\in\ppp(\pr\Gamma)$, we can take $Q\in\pr\Gamma$ such that $\ppp=\hom_\Gamma(Q,\ )$ and $I^0\oplus I^1\in\add(\widehat{Q})^*$ by \XBB. To show that $\ppp$ is minimal, fix $P\in\ind(\pr\Gamma)$. Then $P\in\pp\cc\Leftrightarrow\pd\top P\le 1\Leftrightarrow(\widehat{P})^*\in\add(I^0\oplus I^1)\Rightarrow P\in\add Q\Leftrightarrow\ppp P\in\pr\Lambda$ ($*$) holds. Assume $\ppp P\in\pr\Lambda$, and take $f\in\cc(P^\prime,P)$ such that $\ppp f$ gives the inclusion $J_\Lambda\ppp P\subset \ppp P$. Since $(P]=(0\rightarrow P^\prime\stackrel{f}{\rightarrow}P)$ holds, we obtain $P\in\pp\cc$. Thus above five conditions in ($*$) are equivalent. Consequently, $\ppp$ is minimal by $\add(I^0\oplus I^1)=\add(\widehat{Q})^*$, and we obtain the desired bijections.

(3) Clearly $l^+\supseteq\ind^+_0\cc$ holds. Let $X\in\pp\cc-\ind^+_0\cc$ and $Y:=\nak^+(X)\in\ind\cc$. Then $0\ge l(\mu^+_X)=l(\mu^-_Y)$ holds by \XAD. Since $l(\mu^-_Y)\ge0$ if and only if $Y\in\ii\cc$, we obtain the first equality in (ii). The second equality in (ii) follows from [I3;6.4.1(2)].

Take $I\in\ind(\add(I^0\oplus I^1))$. Let $L:=\soc I$, $0\rightarrow P_1\stackrel{f}{\rightarrow}P_0\rightarrow L\rightarrow0$ a minimal projective resolution and $X:=\ppp P_0\in\ff$. Then $(P_0]=(0\rightarrow P_1\stackrel{f}{\rightarrow}P_0)$ holds for $\pr\Gamma=\cc$. Moreover, $I\in\add I^0$ if and only if $\grade{L}=0$ if and only if $f$ is not an epimorphism in $\cc$ if and only if $P_0\in l^+$ by (ii). Since $\ppp L$ is a simple $\Lambda$-module, we obtain an exact sequence $0\rightarrow\ppp P_1\stackrel{\ppp f}{\longrightarrow}X\rightarrow\top X\rightarrow0$ in $\mod\Lambda$ by taking $\ppp$. Hence $f$ is not an epimorphism in $\cc$ if and only if so is $\ppp f$ in $\ff$ if and only if $\top X\in\ff$. Thus (i) holds.\rule{5pt}{10pt}

\vskip1em{\bf\XCC\ }Now we obtain the following theorem which implies the classical theorem of Auslander in 0.2. Recall that we call an artin algebra $\Gamma$ an {\it Auslander algebra} if $\gl\Gamma\le2$ and $\Gamma$ satisfies the $(1,2)$-condition, namely $\domdim\Gamma\ge2$. Notice that the equivalence of (2) and (3) below is a special case of [R3;Cor. of Prop.6].

\vskip1em{\bf Theorem }{\it Let $\Gamma$ be an artin algebra and $\cc:=\pr\Gamma$. Then the following conditions are equivalent.

(1) $\Gamma$ is an Auslander algebra.

(2) $\Gamma$ is an Auslander-regular ring with $\gl\Gamma\le2$, and any simple $\Gamma$-module $L$ with $\pd L=1$ satisfies $\grade{L}=0$.

(3) There exists an artin algebra $\Lambda$ such that $\cc$ is equivalent to $\mod\Lambda$.

(4) $\cc$ is a strict $\tau$-category and $\nak^-$ gives a map $\ii\cc-\ind^-_0\cc\rightarrow\ind\cc-\pp\cc$.

(5) $\cc$ is a $\tau$-category with an additive function $l$ such that $l^-=\ii\cc$.

(6) $\cc$ is an abelian category.

($i$)$^{op}$ Opposite side version of ($i$) ($1\le i\le 6$).}

\vskip1em{\sc Proof }
Each of above conditions implies that $\Gamma$ is Auslander-regular with $\gl\Gamma\le2$ by \XCA.
Obviously (1) is equivalent to $\add(I^0\oplus I^1)=\add I^0$.
Now \XCB\ immediately implies (1)$\Leftrightarrow$(2)$\Leftrightarrow$(3)$\Leftrightarrow$(4)$^{op}\Leftrightarrow$(5)$^{op}$.
Since $(\mod\Lambda)^{op}$ is equivalent to $\mod\Lambda^{op}$, we obtain (3)$\Leftrightarrow$(3)$^{op}$, and (3)$\Rightarrow$(6) is obvious. We will show (6)$\Rightarrow$(2). We only have to show the latter assertion. Let $0\to P_1\stackrel{f}{\to}P_0\to L\to0$ be a minimal projective resolution. Since $f$ is a non-invertible monomorphism in an abelian category $\cc$, $f$ is not an epimorphism in $\cc$. Thus $\widehat{f}:\widehat{P}_0\to\widehat{P}_1$ is not a monomorphism, and $\widehat{L}\neq0$ holds.\rule{5pt}{10pt}

\vskip1em{\bf\XCCA\ Corollary }{\it
Let $\Lambda$ be a representation-finite artin algebra and $\cc:=\mod\Lambda$. Let $B\in\ind(\pr\Lambda)$, $A:=(\widehat{B})^*\in\ind(\inj\Lambda)$ and $X:=\soc A=\top B\in\ind(\ssim\Lambda)$. Then $(A,\tau^-X)$ is a Nakayama pair if $A$ is not simple, and $(\tau^+X,B)$ is a Nakayama pair if $B$ is not simple.}

\vskip1em{\sc Proof }
Since $\nak^-(A)\in\ind(\mod\Lambda)-\ind(\pr\Lambda)$ holds, there exists an exact sequence $0\rightarrow\tau^+\nak^-(A)\rightarrow A\stackrel{\mu^-_A}{\rightarrow}\theta^-A\rightarrow0$ by the definition \XAD. Since $\mu^-_A$ is a natural surjection $A\rightarrow A/X$, we obtain $\tau^+\nak^-(A)=X$. Thus $(A,\tau^-X)$ is a Nakayama pair.\rule{5pt}{10pt}

\vskip1em{\bf\XCD\ }
We denote by $\mod_{sp}\Lambda$ the full subcategory of $\mod\Lambda$ consisting of $\Lambda$-modules whose socles are projective. Such categories $\mod_{sp}\Lambda$ play an important role in the representation theory. They are closely related to the representation theory of partially ordered sets, vector space categories, and orders over complete discrete vauation rings (see [S]). The theorem below asserts that the endomorphism ring of $\mod_{sp}\Lambda$ is characterized in terms of diagonal Auslander-regular ring, where we call an Auslander-regular ring $\Gamma$ {\it diagonal} if any non-zero direct summand $I$ of $I^i_\Gamma(\Gamma)$ satisfies $\fd_{\Gamma}I=i$.

\vskip1em{\bf Theorem }{\it Let $\Gamma$ be an artin algebra and $\cc:=\pr\Gamma$. Then the following conditions are equivalent.

(1) $\Gamma$ is a diagonal Auslander-regular ring with $\gl\Gamma\le2$.

(2) $\Gamma$ is an Auslander-regular ring with $\gl\Gamma\le2$, and any simple $\Gamma$-module $L$ with $\pd L=1$ satisfies $\grade{L}=1$.

(3) There exists an artin algebra $\Lambda$ such that $\cc$ is equivalent to $\mod_{sp}\Lambda$.

(4) $\cc$ is a strict $\tau$-category and $\nak^-$ gives a (bijective) map $\ii\cc-\ind^-_0\cc\rightarrow\pp\cc-\ind^+_0\cc$.

(5) $\cc$ is a $\tau$-category with an additive function $l$ such that $l^-=\ind^-_0\cc$.

($i$)$^{op}$ Opposite side version of ($i$) ($1\le i\le 5$).}

\vskip1em{\bf\XCDA\ }
Let $\Gamma$ be a $1$-Gorenstein artin algebra.
Then any simple $\Gamma$-module $L$ with $\grade{L}=0$ is projective if and only if any simple $\Gamma^{op}$-module $L$ with $\grade{L}=0$ is projective.

\vskip1em{\sc Proof }
We will show the `only if' part. Assume that a simple $\Gamma^{op}$-module $L$ satisfies $\grade{L}=0$ and $\pd L>0$.
Take a projective resolution $0\rightarrow\Omega L\stackrel{g}{\rightarrow}P_0\stackrel{f}{\rightarrow}L\rightarrow0$.
We have an exact sequence $0\leftarrow\ext_\Gamma^1(L,\Gamma)\stackrel{a}{\leftarrow}\widehat{\Omega}L\stackrel{\widehat{g}}{\leftarrow}\widehat{P}_0\stackrel{\widehat{f}}{\leftarrow}\widehat{L}\leftarrow0$.
Then the injective hull $b\in\hom_\Gamma(\widehat{P}_0,I)$ satisfies $I\in\pr\Gamma$ by the $(1,1)$-condition. Suppose that $b$ is not an isomorphism.
Since $\widehat{b}$ factors through $g$, $b$ factors through $\widehat{g}$. Thus $\widehat{L}=0$, a contradiction. Hence $\widehat{P}_0\in\ind(\inj\Lambda)$ holds. Since $\widehat{\Omega}L\neq0$ by $\pd L>0$, we can take an injection $c\in\hom_\Gamma(M,\widehat{\Omega}L)$ for a simple $\Gamma$-module $M$. Then $M\in\pr\Gamma$ holds by $\grade{M}=0$. Since $ca=0$ holds by the $(1,1)$-condition, there exists $c^\prime$ such that $c=c^\prime\widehat{g}$. Since $\soc\widehat{P}_0$ is simple, $c^\prime$ factors through $\widehat{f}$. Thus $c=0$, a contradiction.\rule{5pt}{10pt}

\vskip1em{\bf\XCDB\ Proof of \XCD\ }
(2)$\Leftrightarrow$(2)$^{op}$ holds by \XCDA. The diagonal condition is equivalent to $\ind(\add(I^0\oplus I^1))-\ind(\add I^1)=\ind(\add I^0)$. Now a similar argument as in the proof of \XCC\ works to show (1)$\Leftrightarrow$(2)$\Leftrightarrow$(3)$\Leftrightarrow$(4)$^{op}\Leftrightarrow$(5)$^{op}$.\rule{5pt}{10pt}

\vskip1em{\bf\XCE\ }We obtain the following corollary by \XABA, where the case $i=2$ is the Auslander's correspondence 0.2.

\vskip1em{\bf Corollary }{\it
There exists a bijection between (1-$i$) and (2-$i$) below ($1\le i\le 4$), which is given by $\cc\mapsto\Gamma:=\cc(M,M)$.

(1) Equivalence classes of additive categories $\cc$ with additive generators $M$ such that

\strut\kern1em(1-1) $\cc$ is a faithful torsionfree class over an artin algebra.

\strut\kern1em(1-2) $\cc$ is a faithful hereditary torsionfree class over an artin algebra.

\strut\kern1em(1-3) $\cc=\mod\Lambda$ over an artin algebra $\Lambda$.

\strut\kern1em(1-4) $\cc=\mod_{sp}\Lambda$ over an artin algebra $\Lambda$.

(2) Morita-equivalence classes of artin algebras $\Gamma$ such that

\strut\kern1em(2-1) $\Gamma$ satisfies $\gl\Gamma\le2$ and the $(2,2)$ and $(2,2)^{op}$-condition.

\strut\kern1em(2-2) $\Gamma$ is an Auslander-regular ring with $\gl\Gamma\le2$.

\strut\kern1em(2-3) $\Gamma$ is an Auslander algebra.

\strut\kern1em(2-4) $\Gamma$ is a diagonal Auslander-regular ring with $\gl\Gamma\le2$.}

\vskip1em{\bf\XCF\ }
In this section, we collect some homological results and questions.

\vskip1em{\bf\XCFA\ Symmetry }
We have obtained a few left-right symmetry in previous sections. Moreover, recall that $\Gamma$ is $n$-Gorenstein if and only if so is $\Gamma^{op}$ by [FGR;3.7], and $\domdim\Gamma=\domdim\Gamma^{op}$ holds by [H2]. These left-right symmetry are generalized as follows: We say that $l\ge0$ is a {\it dominant number} of $\Gamma$ if $\fd I^i_\Gamma(\Gamma)<\fd I^l_\Gamma(\Gamma)$ holds for any $i$ ($0\le i<l$).

\vskip1em{\bf Theorem }{\it
[I9;1.1,2.4] Let $l$ and $n$ be positive integers.

(1) For an $n$-Gorenstein ring $\Gamma$, the set of dominant numbers of $\Gamma$ smaller than $n$ coincides with that of $\Gamma^{op}$. Any dominant number $l$ of $\Gamma$ with $l<n$ satisfies $\fd I^l_\Gamma(\Gamma)=l$. 

(2) Assume that a noetherian ring $\Gamma$ satisfies the $(l,l)$ and $(l,l)^{op}$-condition. Then $\Gamma$ satisfies the $(l,n)$-condition if and only if it satisfies the $(l,n)^{op}$-condition.}

\vskip1em{\bf\XCFB\ Question }Let $\Gamma$ be an artin algebra. Is the condition that $\Gamma$ is diagonal Auslander-regular left-right symmetric? This is true if $\gl\Gamma\le2$ by \XCD. More generally, is the following condition ($*$) left-right symmetric for an $n$-Gorenstein ring $\Gamma$? This is true for $n=1$ by \XCDA.

($*$) Any simple $\Gamma$-module $L$ with $\grade{L}=i$ satisfies $\pd L=i$ for any $i$ ($0\le i<n$).

\vskip1em{\bf\XCFC\ Duality }Let $\Gamma$ be a noetherian ring and $\extf{n}{}:=\ext^n_\Gamma(\ ,\Gamma)$ and $\fff_n:=\soc\extf{n}{}$. Consider the following condition (D$_n$).

\strut\kern1em(D$_n$) $\fff_n$ gives a bijection between isoclasses of simple $\Gamma$-modules $L$ with $\grade{L}=n$ and that of $\Gamma^{op}$. Moreover, $\fff_n\circ\fff_nL$ is isomorphic to $L$, and $\sgrade{\extf{n}{L}/\fff_nL}>n$ holds.

If an artin algebra $\Gamma$ satisfies $\gl\Gamma\le2$ and the $(2,2)$ and $(2,2)^{op}$-condition, then (D$_2$) holds by \XBA\ and \XACB(2). Moreover, if $\Gamma$ is Auslander-regular, then \XCA\ and \XADB(3) imply that (D$_1$) holds, and more strongly $\fff_1$ gives an injection from isoclasses of simple $\Gamma$-modules $L$ with $\grade{L}\le1$ to isoclasses of simple $\Gamma^{op}$-modules.

These observations are generalized as follows:

\vskip1em{\bf Theorem }{\it
Let $n\ge0$ and $\Gamma$ a noetherian algebra satisfying the $(l,l)$ and $(l,l)^{op}$-condition for $l=n,n+1$. Then (D$_n$) holds.
Moreover, if $\gl\Gamma=n\ge2$, then any simple $\Gamma$-module $L$ with $\grade{L}=0$ and $\rgrade{L}=n-1$ satisfies that $\fff_{n-1}L$ is simple and $\sgrade{\extf{n-1}{L}}=n$.}

\vskip1em{\sc Proof }The former assertion was shown in [I9;1.3]. We will show the latter assertion. Put $\cc_{\Gamma}:=\{ M\in\mod\Gamma\ |\ \sgrade{M}\ge n\}$. Then $\cc_{\Gamma}$ (resp. $\cc_{\Gamma^{op}}$) is an abelian subcategory of $\mod\Gamma$ (resp. $\mod\Gamma^{op}$) closed under subfactor modules, and $\extf{n}{}$ gives a duality between $\cc_{\Gamma}$ and $\cc_{\Gamma^{op}}$ such that $\extf{n}{}\circ\extf{n}{}$ is isomorphic to the identity functor by [I5;6.2]. Take a projective resolution $0\rightarrow P_{n-1}\rightarrow\cdots\rightarrow P_0\rightarrow L\rightarrow0$ by $\pd L=n-1$. Taking $\widehat{(\ )}$, we obtain an exact sequence $0\leftarrow \extf{n-1}{L}\leftarrow\widehat{P}_{n-1}\leftarrow\cdots\leftarrow\widehat{P}_0\leftarrow\widehat{L}\leftarrow0$ by $\rgrade{L}=n-1$. Thus $\widehat{L}\in\ind(\pr\Gamma^{op})$ holds by $\gl\Gamma=n$, and we obtain an exact sequence $0\rightarrow L\rightarrow\widehat{\widehat{L}}\rightarrow M\rightarrow0$ with $M:=\extf{n}{\extf{n-1}{L}}\in\cc_\Gamma$ and $\widehat{\widehat{L}}\in\ind(\pr\Gamma)$. Thus $\extf{n-1}{L}=\extf{n}{M}\in\cc_{\Gamma^{op}}$ holds, and $\top M$ is simple. By the remark above, $\fff_{n-1}L=\fff_{n}{M}=\extf{n}{(\top M)}$ is simple.\rule{5pt}{10pt}

\vskip1em{\bf\XCFD\ Question }
When is $\fff_nL$ a simple $\Gamma^{op}$-module for a simple $\Gamma$-module $L$ and $n$?

\vskip1em{\bf\XD\ AR quivers }

In the representation theory, the concept of additive functions often appears. We recall several results below which assert that some representation theoretic diagrams are characterized by the existence of additive functions:

{\small
(a) It is a classical result that Dynkin diagrams and extended Dynkin diagrams are characterized in terms of additive functions [HPR].

(b) Brenner characterized AR quivers of representation-finite artin algebras in terms of {\it hammocks}, which is a formulation of the existence of additive functions [Br]. At the same time, Igusa and Todorov gave another characterization independently which does not use additive functions [IT3].

(c) Ringel and Vossieck studied hammocks in [RVo] very clearly, and characterized AR quivers of representation-finite partially ordered sets in terms of hammocks.

(d) Reiten and Van den Bergh characterized AR quivers of representation-finite two-dimensional orders, which essentially uses additive functions [RV].

(e) Inspired by the work of Igusa-Todorov (b), the author characterized AR quivers of representation-finite one-dimensional orders [I5], which does not use additive functions. Then Rump gave another characterization in terms of additive functions [R5].

(f) Additive functions are used to characterize rejectable subsets \XEA(2) for two-dimensional orders by Reiten and Van den Bergh [RV] and for one-dimensional orders by the author [I2].}

In this section, we shall see that (b) and (c) above are understood clearly in our viewpoint of (e) and this paper (see \XDDA).

\vskip1em{\bf\XDA\ Definition }
(1) $Q=(Q,Q^p,Q^i,\tau^+,d,d^\prime)$ is called a {\it translation quiver} if $Q$ is a set, $Q^p$ and $Q^i$ are subsets of $Q$, $\tau^+$ is a bijection $Q-Q^p\rightarrow Q-Q^i$, $d$ and $d^\prime$ are maps $Q\times Q\rightarrow\nnn_{\ge0}$ such that $d(Y,X)=d^\prime(\tau^+X,Y)$ holds for any $X\in Q-Q^p$ and $Y\in Q$, and $d(\ ,X)=0$ implies $X\in Q^p$. We call $Q$ {\it admissible} if there exists a map $c:Q\rightarrow\nnn_{>0}$ such that $c(X)d(X,Y)=d^\prime(X,Y)c(Y)$ holds for any $X$, $Y\in Q$. We call $Q$ {\it locally finite} if $\sum_{Y\in Q}d(Y,X)<\infty$ and $\sum_{Y\in Q}d^\prime(X,Y)<\infty$ hold for any $Y\in Q$.

Usually, we draw $Q$ as a directed graph:
$Q$ is the set of vertices, and we draw valued arrows $X\stackrel{(d(X,Y),d^\prime(X,Y))}{\begin{picture}(60,8)
\put(0,4){\vector(1,0){60}}
\end{picture}}Y$
for any $X,Y\in Q$ such that $d(X,Y)\neq0$, and dotted arrows from $X$ to $\tau^+X$ for any $X\in Q-Q^p$.

For a $\tau$-category $\cc$, we define a locally finite translation quiver $\AR(\cc)=(Q,Q^p,Q^i,\tau^+,d,d^\prime)$ called the {\it AR quiver} of $\cc$ as follows: $Q:=\ind\cc$, $Q^p:=\pp\cc$, $Q^i:=\ii\cc$, $d(X,Y)$ is the multiplicity of $X$ in $\theta^+Y$ and $d^\prime(X,Y)$ is the multiplicity of $Y$ in $\theta^-X$. Thus $\AR(\cc)$ displays terms of each $(X]$ and $[X)$ diagrammatically. If $\cc$ is a torsionfree class over an artin algebra $\Lambda$ over $R$, then $\AR(\cc)$ is admissible by $k:=R/J_R$ and $c(X):=\dim_k\endm_\Lambda(X)/J_{\endm_\Lambda(X)}$.

(2) For a locally finite translation quiver $Q$, we denote by $\zzz Q$ (resp. $\nnn Q$) the free $\zzz$-module (resp. free abelian monoid) generated by $Q$. For $X=\sum_{Y\in Q}a_YY\in\zzz Q$, put $\supp Y:=\{ X\in Q\ |\ a_X\neq0\}$. Define elements $\theta^+$, $\theta^-$, $\tau^+$ and $\tau^-$ of $\endm_{\zzz}(\zzz Q)$ as follows: Put $\theta^+X:=\sum_{Y\in Q}d(Y,X)Y$ and $\theta^-X:=\sum_{Y\in Q}d^\prime(X,Y)Y$ for $X\in Q$. Put $\tau^+X:=0$ for $X\in Q^p$, $\tau^-X:=(\tau^+)^{-1}(X)$ for $X\in Q-Q^i$ and $\tau^-X:=0$ for $X\in Q^i$. When $Q=\AR(\cc)$ for a $\tau$-category $\cc$, these definitions are consistent with those in \XACA. Define $\phi^\pm$ and a ({\it left, right}) {\it additive function} of $Q$ by a similar manner in \XACA, and define $\theta^+_n$ and $\eta^+_n\in\endm_{\zzz}(\zzz Q)$ ($n\ge0$) by the recursion formulas in \XACC\ and \XAD.

\vskip1em{\bf\XDB\ }We call an additive category $\cc$ with an additive generator $M$ {\it artinian} if the ring $\cc(M,M)$ is artinian.
We call a translation quiver $Q$ {\it artinian} (resp. {\it strict}) if there exists an artinian (resp. strict) $\tau$-category $\cc$ with $Q=\AR(\cc)$.
The proposition below gives a simple criterion for $Q$ to be artinian (resp. strict).

\vskip1em{\bf Proposition }{\it
Let $Q$ be an admissible translation quiver with a finite number of vertices. Then $Q$ is artinian (resp. strict) if and only if any $\tau$-category $\cc$ with $Q=\AR(\cc)$ is artinian (resp. strict).

(1) $Q$ is artinian if and only if there exists $n>0$ such that $\theta^+_n=0$.
If any connected component of $Q$ contains a vertex in $Q^i$, then $Q$ is artinian if and only if there exists $n>0$ such that $\theta^+_nX=0$ for any $X\in Q^i$.

(2) If $Q$ is artinian, then $Q$ is strict if and only if $Q=\bigcup_{X\in Q^i,n\ge0}\supp\theta^+_nX$ if and only if $\theta^+_n=\theta^+\circ\theta^+_{n-1}-\tau^+\circ\theta^+_{n-2}$ for any $n\ge2$.}

\vskip1em{\bf\XDBA\ Example }Let $Q$ be the translation quiver below, where $\tau^+$ is the left translation.
{\scriptsize\[\begin{array}{ccccccccccccc}
\Q{1}&&&&\Q{13}&&&&\Q{6}&&&&\Q{1}\\
&\QDR&&\QUR&&\QDR&&\QUR&&\QDR&&\QUR&\\
&&\Q{18}&\QR&\Q{14}&\QR&\Q{11}&\QR&\Q{7}&\QR&\Q{4}&&\\
&\QUR&&\QDR&&\QUR&&\QDR&&\QUR&&\QDR&\\
\Q{2}&&&&\Q{15}&&&&\Q{8}&&&&\Q{2}\\
&\QDR&&\QUR&&\QDR&&\QUR&&\QDR&&\QUR&\\
&&\Q{19}&\QR&\Q{16}&\QR&\Q{12}&\QR&\Q{9}&\QR&\Q{5}&&\\
&\QUR&&\QDR&&\QUR&&\QDR&&\QUR&&\QDR&\\
\Q{3}&&&&\Q{17}&&&&\Q{10}&&&&\Q{3}
\end{array}\ \ \ \ \ \begin{array}{l}
Q^p=\{14,16\}\\
Q^i=\{7,9\}
\end{array}\]}

Then the calculation of $\theta^+_n(7)$ and $\theta^+_n(9)$ below implies that $Q$ is artinian by \XDB(1), where we describe the diagram in \XACC.
{\scriptsize\begin{eqnarray*}
&\begin{array}{ccccccccccccccc}
16&\rightarrow&12&\rightarrow&8,10&\rightarrow&4,5&\rightarrow&1,2&\rightarrow&18&\rightarrow&14&\rightarrow&\\
\downarrow&&\downarrow&&\downarrow&&\downarrow&&\downarrow&&\downarrow&&\downarrow&&\downarrow\\
&\rightarrow&9&\rightarrow&5&\rightarrow&2,3&\rightarrow&18,19&\rightarrow&13,15&\rightarrow&11&\rightarrow&7\end{array}&\\
&\begin{array}{ccccccccccccccc}
14&\rightarrow&11&\rightarrow&6,8&\rightarrow&4,5&\rightarrow&2,3&\rightarrow&19&\rightarrow&16&\rightarrow&\\
\downarrow&&\downarrow&&\downarrow&&\downarrow&&\downarrow&&\downarrow&&\downarrow&&\downarrow\\
&\rightarrow&7&\rightarrow&4&\rightarrow&1,2&\rightarrow&18,19&\rightarrow&15,17&\rightarrow&12&\rightarrow&9\end{array}&\end{eqnarray*}}

\vskip-1em
On the other hand, $Q$ is not strict by \XDB(2) and $\bigcup_{n\ge0}\supp\theta^+_n(7)\cup\supp\theta^+_n(9)=Q-\{6,8,10,14,16\}$.

\vskip1em{\bf\XDBB\ Lemma }{\it Let $\cc$ be a $\tau$-category and $Q=\AR(\cc)$. Assume that any connected component of $Q$ contains a vertex in $Q^i$. Then $\cc(\ ,X)$ has finite length for any $X\in\cc$ if and only if $\cc(\ ,X)$ has finite length for any $X\in\ii\cc$.}

\vskip 1em{\sc Proof }
Put $\dd:=\{ X\in\cc|\cc(\ ,X)\mbox{ has finite length}\}$, $\underline{\cc}:=\cc/[\pp\cc]$, $\overline{\cc}:=\cc/[\ii\cc]$ and $\cc^\prime:=\add(\ii\cc)$. Then $\cc^\prime\subseteq\dd$ holds by our assumption. Fix indecomposable $X\in\dd$. Since $\underline{\cc}(\ ,X)$ has finite length by $X\in\dd$, $\overline{\cc}(\ ,\tau^+X)$ has finite length by the proof of [I4;2.4]. Now we will show $Y:=\tau^+X\in\dd$. Since we have an exact sequence $0\rightarrow[\cc^\prime](\ ,Y)\rightarrow\cc(\ ,Y)\rightarrow\overline{\cc}(\ ,Y)\longrightarrow0$, we only have to show that $[\cc^\prime](\ ,Y)$ has finite length. Since any finite length $\cc$-module is finitely presented, $[\cc^\prime](\ ,Y)$ is a finitely generated $\cc$-module. Thus there exists $f\in\cc(Z,Y)$ with $Z\in\cc^\prime$ such that $\cc(\ ,Z)\stackrel{\cdot f}{\rightarrow}[\cc^\prime](\ ,Y)\rightarrow0$ is exact. Hence $[\cc^\prime](\ ,Y)$ has finite length by $\cc^\prime\subseteq\dd$. Thus $\tau^+X\in\dd$ holds. Since we have an exact sequence $\cc(\ ,\tau^+X)\rightarrow\cc(\ ,\theta^+X)\rightarrow\cc(\ ,X)$, we obtain $\theta^+X\in\dd$. Thus any predecessor of $X$ in $\AR(\cc)$ is again contained in $\dd$. By our assumption, we can easily show that there exists a path from any vertex in $Q$ to some vertex in $Q^i$. Thus $\dd=\cc$ holds.\rule{5pt}{10pt}

\vskip1em{\bf\XDBC\ Proof of \XDB\ }
(1)(cf. [I3;7.3]) The former assertion follows from \XACC, and the latter assertion follows from the former one and \XDBB.\ \ \ (2) follows from [I3;7.4].

Now the first equivalence follows from (1) and (2), where we remark that any admissible translation quiver is realized as the AR quiver of some $\tau$-category [I5;4.2.1].\rule{5pt}{10pt}

\vskip1em{\bf\XDC\ Theorem }{\it
Let $\Gamma$ be an artin algebra and $\cc:=\pr\Gamma$.
Put $l_X(Y):=\length_{\cc(X,X)}\cc(X,Y)$ for any $X,Y\in\cc$.

(1) Assume that the conditions in \XBA\ are satisfied. Then $\phi^+(\nnn\ind\cc)\supseteq\nnn\ind\cc$ holds. Moreover, a map $l:\ind\cc\rightarrow\nnn_{>0}$ is a right additive function if and only if $l=\sum_{X\in\pp\cc}a_Xl_X$ holds for some $(a_X)\in\nnn^{\pp\cc}$. Such $(a_X)$ is uniquely determined.

(2) Assume that the conditions in \XCA\ are satisfied and put $\dn{S}^+(\cc):=\{ X\in\pp\cc\ |\ \mu^+_X$ is not an epimorphism$\}$. Then a map $l:\ind\cc\rightarrow\nnn_{>0}$ is an additive function if and only if $l=\sum_{X\in\dn{S}^+(\cc)}a_Xl_X$ holds for some $(a_X)\in\nnn^{\dn{S}^+(\cc)}$. Such $(a_X)$ is uniquely determined.}

\vskip1em{\sc Proof }
Since $0\rightarrow\cc(X,\tau^+Y)\rightarrow\cc(X,\theta^+Y)\rightarrow\jac_{\cc}(X,Y)\rightarrow0$ is exact for any $Y\in\cc$, we obtain $l_X(\phi^+Y)=0$ for any $Y\in\ind\cc-\{X\}$ and $l_X(\phi^+X)=1$.

(1) Since $\cc$ is artinian strict, $\phi^+(\sum_{n\ge0}\theta^+_nX)=\sum_{n\ge0}\theta^+_nX-\sum_{n\ge0}(\theta^+\theta^+_nX-\tau^+\theta^+_{n-1}X)=\sum_{n\ge0}\theta^+_nX-\sum_{n\ge0}\theta^+_{n+1}X=X$ holds for any $X\in\cc$ by \XDB(2). Thus the first assertion follows. We will show the `only if' part. Put $a_X:=l(\phi^+X)\ge0$ for any $X\in\pp\cc$ and $l^\prime:=l-\sum_{X\in\pp\cc}a_Xl_X$. Since $l^\prime\circ\phi^+=0$ holds, we obtain $l^\prime=0$ by the first assertion. Thus $l=\sum_{X\in\pp\cc}a_Xl_X$ holds. The third assertion follows from the first.

(2) Since $l(\phi^+X)=0$ holds for any $X\in\pp\cc-\dn{S}^+(\cc)$ by \XCB(3), we obtain $l=\sum_{X\in\dn{S}^+(\cc)}a_Xl_X$ by (1). Fix $X\in\dn{S}^+(\cc)$. We only have to show that $l_X(\mu^-_Y)\ge0$ holds for any $Y\in\ii\cc$.
Put $Z:=\nak^-(Y)\in\ind\cc$, then $l_X(\mu^-_Y)=l_X(\mu^+_Z)$ holds by \XAD. Since $\nak^+(X)\notin\ii\cc$ holds by \XCB(3), we obtain $X\neq Z$. Thus $l_X(\phi^+Z)=0$ implies $l_X(\mu^+_Z)=l_X(\tau^+Z)\ge0$.\rule{5pt}{10pt}

\vskip1em{\bf\XDD\ Theorem }{\it
Let $Q$ be an admissible artinian translation quiver with a finite number of vertices. Then the following conditions (1-i), (2-i), (3-i) and (4-i) are equivalent for each $i$ ($1\le i\le 4$).

(1) There is an artin algebra $\Lambda$

\strut\kern1em
(1-1) with a torsion theory $(\tt,\ff)$ on $\mod\Lambda$ such that $Q=\AR(\ff)$.

\strut\kern1em
(1-2) with a hereditary torsion theory $(\tt,\ff)$ on $\mod\Lambda$ such that $Q=\AR(\ff)$.

\strut\kern1em
(1-3) such that $Q=\AR(\mod\Lambda)$.

\strut\kern1em
(1-4) such that $Q=\AR(\mod_{sp}\Lambda)$.

(2) Any $\tau$-category $\cc$ with $\AR(\cc)=Q$ is

\strut\kern1em
(2-1) strict.

\strut\kern1em
(2-2) strict and $\nak^-$ gives a map $\ii\cc-\ind^-_0\cc\rightarrow\ind\cc$.

\strut\kern1em
(2-3) strict and $\nak^-$ gives a map $\ii\cc-\ind^-_0\cc\rightarrow\ind\cc-\pp\cc$.

\strut\kern1em
(2-4) strict and $\nak^-$ gives a map $\ii\cc-\ind^-_0\cc\rightarrow\pp\cc-\ind^+_0\cc$.

(3)(3-1) $Q=\bigcup_{X\in Q^i,n\ge0}\supp\theta^+_nX$.

\strut\kern1em
(3-2) $Q=\bigcup_{X\in Q^i,n\ge0}\supp\theta^+_nX$, and for any $X\in Q^i$ with $\theta^-X\neq0$, there exists $n\ge0$ such that $\eta^+_iX\in\nnn(Q-Q^p)$ for any $i$ ($0\le i<n$) and $\eta^+_{n+1}X=0$.

\strut\kern1em
(3-3) (3-2) holds, and $\eta^+_nX\in Q-Q^p$.

\strut\kern1em
(3-4) (3-2) holds, and $\eta^+_nX\in Q^p$.

(4)(4-1) $Q$ has a right additive function $l$.

\strut\kern1em
(4-2) $Q$ has an additive function $l$.

\strut\kern1em
(4-3) $Q$ has an additive function $l$ with $l^-=\ii\cc$.

\strut\kern1em
(4-4) $Q$ has an additive function $l$ with $l^-=\ind^-_0\cc$.}

\vskip1em{\bf\XDDA\ Remark }(1) In \XDD(2), we can replace `any' by `some'. In \XDD(3-2), $\nak^-(X)=\eta^+_nX$ holds. In \XDD(1), we can add the condition that $\Lambda$ is a finite dimensional algebra over arbitrary finite field $k$ (see [I5;4.2.1]).

(2) Our condition (3-3) simplifies that of Igusa-Todorov (b) above, and our condition (4-3) simplifies that of Brenner (b) above. After the work of Brenner, Ringel and Vossieck (c) call a simply connected translation quiver $Q$ with a unique source $X$ {\it hammock} if $Q$ has an additive function $l$ with $l^-=\{ X\}$. Thus our condition (4-4) is a generalization of their hammock condition to a general translation quivers.

\vskip1em{\bf\XDDB\ Proof of \XDD\ }
We can fix an artin algebra $\Gamma_0$ such that $\cc_0:=\pr\Gamma_0$ forms a $\tau$-category with $Q=\AR(\cc_0)$ by [I5;4.2.1]. Then (1-$i$)$\Rightarrow$(4-$i$)$\Rightarrow$(2-$i$)$\Rightarrow$($\cc:=\cc_0$ satisfies (2-$i$))$\Rightarrow$(1-$i$) holds by \XBA\ ($i=1$), \XCA\ ($i=2$), \XCC\ ($i=3$) and \XCD\ ($i=4$). Moreover, (2-$i$)$\Leftrightarrow$(3-$i$) holds by \XDB(1) ($i=1$) and \XADB(1) ($i=2,3,4$).\rule{5pt}{10pt}


\vskip1em{\bf\XDE\ Example }
(1) Let $Q$ be the artinian strict translation quiver below, where $\tau^+$ is the left translation.

{\scriptsize\[\begin{array}{cccccccccccccccccccccccccc}
\Q{1}&&&&\Q{33}&&&&\Q{27}&\QR&\Q{24}&\QR&\Q{20}&&&&\Q{14}&&&&\Q{9}&\QR&\Q{4}&\QR&\Q{1}\\
&\QDR&&\QUR&&\QDR&&\QUR&&\QDR&&\QUR&&\QDR&&\QUR&&\QDR&&\QUR&&\QDR&&\QUR&\\
&&\Q{37}&\QR&\Q{34}&\QR&\Q{31}&\QR&\Q{28}&\QR&\Q{25}&\QR&\Q{21}&\QR&\Q{18}&\QR&\Q{15}&\QR&\Q{12}&&&&\Q{5}&&\\
&\QUR&&\QDR&&\QUR&&\QDR&&\QUR&&\QDR&&\QUR&&\QDR&&\QUR&&\QDR&&\QUR&&\QDR&\\
\Q{2}&&&&\Q{35}&&&&\Q{29}&&&&\Q{22}&&&&\Q{16}&&&&\Q{10}&\QR&\Q{6}&\QR&\Q{2}\\
&\QDR&&\QUR&&\QDR&&\QUR&&\QDR&&\QUR&&\QDR&&\QUR&&\QDR&&\QUR&&\QDR&&\QUR&\\
&&\Q{38}&&&&\Q{32}&&&&\Q{26}&&&&\Q{19}&&&&\Q{13}&&&&\Q{7}&&\\
&\QUR&&\QDR&&\QUR&&\QDR&&\QUR&&\QDR&&\QUR&&\QDR&&\QUR&&\QDR&&\QUR&&\QDR&\\
\Q{3}&&&&\Q{36}&&&&\Q{30}&&&&\Q{23}&&&&\Q{17}&&&&\Q{11}&\QR&\Q{8}&\QR&\Q{3}
\end{array}\ \ \ \ \ \begin{array}{c}Q^p=\{4,6,8,24,34\}\\
Q^i=\{4,6,8,15,24\}
\end{array}\]}

Then the calculation of $\eta^+_n(Q^i)$ below shows that $Q$ satisfies the condition \XDD(3-2), where we describe the diagram in \XAD. Thus $Q=\AR(\cc)$ holds for some hereditary torsionfree class $\cc$ over an artin algebra $\Lambda$.
{\scriptsize\begin{eqnarray*}&\begin{array}{ccccccccc}
19&\rightarrow&16
&\rightarrow&12
&\rightarrow&9
&\rightarrow&4\\
\downarrow&&\downarrow&&\downarrow&&\downarrow&&\downarrow\\
17&\rightarrow&13
&\rightarrow&10
&\rightarrow&5
&\rightarrow&1
\end{array}&\nak^-(4)=17\\
&\begin{array}{ccccccccccccccc}
27
&\rightarrow&25
&\rightarrow&21,22
&\rightarrow&18,19
&\rightarrow&14,16,17
&\rightarrow&12,13
&\rightarrow&10
&\rightarrow&6\\
\downarrow&&\downarrow&&\downarrow&&\downarrow&&\downarrow&&\downarrow&&\downarrow&&\downarrow\\
24
&\rightarrow&20
&\rightarrow&18
&\rightarrow&15,16
&\rightarrow&12,13
&\rightarrow&9,10,11
&\rightarrow&5,7
&\rightarrow&2
\end{array}&\nak^-(6)=24\\
&\begin{array}{ccccccccccccccccc}
32&\rightarrow&29
&\rightarrow&25
&\rightarrow&20,21
&\rightarrow&18
&\rightarrow&16
&\rightarrow&13
&\rightarrow&11
&\rightarrow&8\\
\downarrow&&\downarrow&&\downarrow&&\downarrow&&\downarrow&&\downarrow&&\downarrow&&\downarrow&&\downarrow\\
30&\rightarrow&26
&\rightarrow&22
&\rightarrow&18
&\rightarrow&14,15
&\rightarrow&12
&\rightarrow&10
&\rightarrow&7
&\rightarrow&3
\end{array}&\nak^-(8)=30\\
&\begin{array}{ccccccccccccccccccc}
37
&\rightarrow&33,35
&\rightarrow&31,32
&\rightarrow&28,29,30
&\rightarrow&25,26
&\rightarrow&20,22
&\rightarrow&18
&\rightarrow&15\\
\downarrow&&\downarrow&&\downarrow&&\downarrow&&\downarrow&&\downarrow&&\downarrow&&\downarrow\\
34
&\rightarrow&31
&\rightarrow&27,29
&\rightarrow&25,26
&\rightarrow&21,22,23
&\rightarrow&18,19
&\rightarrow&14,16
&\rightarrow&12
\end{array}&\nak^-(15)=34\\
&\begin{array}{ccccccccccccccc}
10
&\rightarrow&5,7
&\rightarrow&1,2,3
&\rightarrow&37,38
&\rightarrow&34,35
&\rightarrow&31
&\rightarrow&27
&\rightarrow&24\\
\downarrow&&\downarrow&&\downarrow&&\downarrow&&\downarrow&&\downarrow&&\downarrow&&\downarrow\\
6
&\rightarrow&2
&\rightarrow&37,38
&\rightarrow&33,35,36
&\rightarrow&31,32
&\rightarrow&28,29
&\rightarrow&25
&\rightarrow&20
\end{array}&\nak^-(24)=6\end{eqnarray*}}

(2) Let $Q$ be the artinian strict translation quiver below, where $\tau^+$ is the left translation.
{\scriptsize\[\begin{array}{cccccccccccccccccc}
&&&&&&&&\Q{13}&&&&\Q{7}&&&&\Q{1}\\
&&&&&&&&&\QDR&&\QUR&&\QDR&&\QUR&\\
&&&&&&&&&&\Q{11}&&&&\Q{5}&\\
&&&&&&&&&\QUR&&\QDR&&\QUR&&\QDR&\\
\Q{3}&&&&\Q{20}&&&&\Q{14}&&&&\Q{8}&&&&\Q{2}\\
&\QDR&&\QUR&&\QDR&&\QUR&&\QDR&&\QUR&&\QDR&&\QUR&\\
\Q{4}&\QR&\Q{24}&\QR&\Q{21}&\QR&\Q{18}&\QR&\Q{15}&\QR&\Q{12}&\QR&\Q{9}&\QR&\Q{6}&\QR&\Q{3}\\
&\QUR&&\QDR&&\QUR&&\QDR&&\QUR&&\QDR&&\QUR&&\QDR&\\
\Q{2}&&&&\Q{22}&&&&\Q{16}&&&&\Q{10}&&&&\Q{4}\\
&\QDR&&\QUR&&\QDR&&\QUR&&&&&&&&&\\
&&\Q{25}&&&&\Q{19}&&&&&&&&&&\\
&\QUR&&\QDR&&\QUR&&\QDR&&&&&&&&&\\
\Q{26}&&&&\Q{23}&&&&\Q{17}&&&&&&&&
\end{array}\ \ \ \ \ \begin{array}{c}Q^p=\{11,13,25,26\}\\
Q^i=\{1,5,17,19\}\end{array}\]}

Then the calculation of $\eta^+_n(Q^i)$ below shows that $Q$ satisfies the condition \XDD(3-4), where we describe the diagram in \XAD. Thus $Q=\AR(\mod_{sp}\Lambda)$ holds for some artin algebra $\Lambda$.
{\scriptsize\begin{eqnarray*}
&\begin{array}{ccccccccccccccc}
2,26
&\rightarrow&24,25
&\rightarrow&20,21,22
&\rightarrow&18^2
&\rightarrow&14,15,16
&\rightarrow&11,12
&\rightarrow&7,8
&\rightarrow&5\\
\downarrow&&\downarrow&&\downarrow&&\downarrow&&\downarrow&&\downarrow&&\downarrow&&\downarrow\\
25
&\rightarrow&22,23
&\rightarrow&18,19
&\rightarrow&14,15,16
&\rightarrow&12^2
&\rightarrow&8,9,10
&\rightarrow&5,6
&\rightarrow&1,2
\end{array}&\nak^-(5)=25\\
&\begin{array}{ccccccccccccccc}
13,14
&\rightarrow&11,12
&\rightarrow&8,9,10
&\rightarrow&6^2
&\rightarrow&2,3,4
&\rightarrow&24,25
&\rightarrow&22,23
&\rightarrow&19\\
\downarrow&&\downarrow&&\downarrow&&\downarrow&&\downarrow&&\downarrow&&\downarrow&&\downarrow\\
11
&\rightarrow&7,8
&\rightarrow&5,6
&\rightarrow&2,3,4
&\rightarrow&24^2
&\rightarrow&20,21,22
&\rightarrow&18,19
&\rightarrow&16,17
\end{array}&\nak^-(19)=11\end{eqnarray*}}

(3) Let $Q$ be the artinian strict translation quiver below, where $\tau^+$ is the left translation.
{\scriptsize\[\begin{array}{ccccccccccccc}
\Q{1}&\QR&\Q{17}&&&&\Q{11}&&&&\Q{5}&\QR&\Q{1}\\
&\QUR&&\QDR&&\QUR&&\QDR&&\QUR&&\QDR&\\
\Q{2}&\QR&\Q{18}&\QR&\Q{15}&\QR&\Q{12}&\QR&\Q{9}&\QR&\Q{6}&\QR&\Q{2}\\
&\QDR&&\QUR&&\QDR&&\QUR&&\QDR&&\QUR&\\
&&\Q{19}&&&&\Q{13}&&&&\Q{7}&&\\
&\QUR&&\QDR&&\QUR&&\QDR&&\QUR&&\QDR&\\
\Q{3}&&&&\Q{16}&&&&\Q{10}&&&&\Q{3}\\
&\QDR&&\QUR&&\QDR&&\QUR&&\QDR&&\QUR&\\
\Q{4}&\QR&\Q{20}&&&&\Q{14}&&&&\Q{8}&\QR&\Q{4}
\end{array}\ \ \ \ \ Q^p=\{1,4\}=Q^i\]}

Then the calculation of $\eta^+_n(Q^i)$ below shows that $Q$ satisfies the condition \XDD(3-3), where we describe the diagram in \XAD.
{\scriptsize\begin{eqnarray*}
&\begin{array}{ccccccccccccccccc}
9&\rightarrow&6,  7&\rightarrow&2,  3&\rightarrow&17,  19,  20&\rightarrow&15,  16&\rightarrow&12,  13&\rightarrow&9&\rightarrow&5&\rightarrow&1\\
\downarrow&&\downarrow&&\downarrow&&\downarrow&&\downarrow&&\downarrow&&\downarrow&&\downarrow&&\downarrow\\
5&\rightarrow&2&\rightarrow&18,  19&\rightarrow&15,  16&\rightarrow&11,  13,  14&\rightarrow&9,  10&\rightarrow&6,  7&\rightarrow&2&\rightarrow&17
\end{array}&\nak^-(1)=5\\
&\begin{array}{ccccccccccccccccc}
10&\rightarrow&7&\rightarrow&2&\rightarrow&17,  18&\rightarrow&15&\rightarrow&13&\rightarrow&10&\rightarrow&8&\rightarrow&4\\
\downarrow&&\downarrow&&\downarrow&&\downarrow&&\downarrow&&\downarrow&&\downarrow&&\downarrow&&\downarrow\\
8&\rightarrow&3&\rightarrow&19&\rightarrow&15&\rightarrow&11,  12&\rightarrow&9&\rightarrow&7&\rightarrow&3&\rightarrow&20
\end{array}&\nak^-(4)=8
\end{eqnarray*}}

For example, $Q=\AR(\mod\Lambda)$ holds for a non-standard algebra $\Lambda$ which is defined by the quiver {\scriptsize$\stackrel{b}{\begin{picture}(20,10)
\put(20,0){\line(-1,0){15}}
\put(5,3){\oval(6,6)[l]}
\put(5,6){\vector(1,0){15}}
\end{picture}}\bullet\begin{array}{c}\stackrel{a}{\longrightarrow}\\
\stackrel{c}{\longleftarrow}\end{array}\bullet$} with relations $\{ ca=b^2,b^4=0,ac=abc\}$ [BG].

(4) Let $Q$ be the artinian strict translation quiver below, where $\tau^+$ is the left translation. Then $Q$ is a part of the preprojective component of the hereditary algebra $\Lambda$ of the wild quiver {\scriptsize$\begin{array}{ccl}
&\bullet&\bullet\\
&\downarrow&\uparrow\\
\bullet&\rightarrow\bullet\leftarrow&\bullet\rightarrow\bullet\leftarrow\bullet\end{array}$}, and $Q=\AR(\cc)$ holds for a faithful torsionfree class $\cc$ of $\mod\Lambda$.
{\scriptsize\[\begin{array}{cccccccccccccccccc}
\Q{22}&&&&\Q{15}&&&&\Q{8}&&&&&&&&\\
&\QDR&&\QUR&&\QDR&&\QUR&&\QDR&&&&&&&\\
\Q{23}&\QR&\Q{19}&\QR&\Q{16}&\QR&\Q{12}&\QR&\Q{9}&\QR&\Q{5}&&&&&&\\
&\QUR&&\QDR&&\QUR&&\QDR&&\QUR&&\QDR&&&&&\\
\Q{24}&\QR&\Q{20}&\QR&\Q{17}&\QR&\Q{13}&\QR&\Q{10}&\QR&\Q{6}&\QR&\Q{3}&&&&\\
&\QDR&&\QUR&&\QDR&&\QUR&&\QDR&&\QUR&&\QDR&&&\\
&&\Q{21}&&&&\Q{14}&&&&\Q{7}&&&&\Q{2}&&\\
&\QUR&&\QDR&&\QUR&&\QDR&&\QUR&&\QDR&&\QUR&&\QDR&\\
\Q{25}&&&&\Q{18}&&&&\Q{11}&&&&\Q{4}&&&&\Q{1}
\end{array}\ \ \ \ \ \begin{array}{c}Q^p=\{19,20,21,22,23,24,25\}\\
Q^i=\{1,2,3,5,6,8,9\}\end{array}\]}

In fact, $Q$ satisfies the condition \XDD(3-1), bot not (3-2) since the calculation below shows that $\nak^-(i)$ is not defined for $i=5,6,8,9$.
{\scriptsize\begin{eqnarray*}
&\begin{array}{ccccccccccccc}
&\rightarrow&22,23,24
&\rightarrow&19,21
&\rightarrow&17,18
&\rightarrow&13,14
&\rightarrow&10
&\rightarrow&5\\
\downarrow&&\downarrow&&\downarrow&&\downarrow&&\downarrow&&\downarrow&&\downarrow\\
22,23,24^2
&\rightarrow&19^2,20
&\rightarrow&15,16,17
&\rightarrow&12,14
&\rightarrow&10,11
&\rightarrow&6,7
&\rightarrow&3
\end{array}&\\
&\begin{array}{ccccccccccccc}
&\rightarrow&22,23,24^3
&\rightarrow&19^2,20,21
&\rightarrow&15,16,17,18
&\rightarrow&12,14
&\rightarrow&10
&\rightarrow&6\\
\downarrow&&\downarrow&&\downarrow&&\downarrow&&\downarrow&&\downarrow&&\downarrow\\
22^2,23^2,24^4,25^2
&\rightarrow&19^3,20^2,21^2
&\rightarrow&15,16,17^3
&\rightarrow&12^2,13,14
&\rightarrow&8,9,10,11
&\rightarrow&5,7
&\rightarrow&3
\end{array}&\\
&\begin{array}{ccccccccccc}
&\rightarrow&22,24^2,25
&\rightarrow&19,20,21
&\rightarrow&16,17
&\rightarrow&12
&\rightarrow&8\\
\downarrow&&\downarrow&&\downarrow&&\downarrow&&\downarrow&&\downarrow\\
22,23^2,24^3,25
&\rightarrow&19^2,20,21^2
&\rightarrow&15,17^2,18
&\rightarrow&12,13,14
&\rightarrow&9,10
&\rightarrow&5
\end{array}&\\
&\begin{array}{ccccccccccc}
&\rightarrow&23,24^2,25
&\rightarrow&19,20,21
&\rightarrow&15,17
&\rightarrow&12
&\rightarrow&9\\
\downarrow&&\downarrow&&\downarrow&&\downarrow&&\downarrow&&\downarrow\\
22^2,23,24^3,25
&\rightarrow&19^2,20,21^2
&\rightarrow&16,17^2,18
&\rightarrow&12,13,14
&\rightarrow&8,10
&\rightarrow&5
\end{array}&
\end{eqnarray*}}

\vskip0em{\bf\XE\ Rejection theory}

The {\it rejection theory} of an additive category $\cc$ is a study of subcategories $\cc^\prime$ called {\it rejective} (\XEA) and the corresponding subset $\ind\cc^\prime$ of $\ind\cc$ called {\it rejectable} (\XEA). For example, rejective subcategories of $\mod\Lambda$ for an artin algebra $\Lambda$ are given by factor algebras of $\Lambda$ (\XEAA). The first example of rejection theory seems to be DK (=Drozd-Kirichenko) Rejection Lemma [DK] (see \XECA), which characterizes one-point rejectable subsets and plays a crucial role in the theory of Bass orders [DKR][Ro][HN]. In [I1], the author studied the rejection theory of orders and artin algebras by connecting with Auslander-Reiten theory, and characterized finite rejectable subsets in terms of AR quivers (see \XEC), a generalization of DK Rejection Lemma.
It is surprizing that the rejectability of a finite subset ${\bf S}$ of $\ind\cc$ depends only on the restriction of $\AR(\cc)$ to ${\bf S}$ (see [I1] for examples of rejectable subsets). Moreover, he studied the rejection theory of arbitrary $\tau$-categories in [I4]. In particular, he generalized results in [I1] to $\tau$-categories, and successfully applied to characterize AR quivers of representation-finite orders [I5].
Of course, his results are valid for our case when $\cc$ is a torsionfree class of $\mod\Lambda$, and we will give a representation theoretic interpretation of rejective subcategories of such $\cc$ in \XEB.
Note that, recently rejective subcategories was successfully applied to quite different kind of problems, Solomon's second conjecture on zeta functions of orders [I7] and finiteness of representation dimension of artin algebras [I8].

\vskip1em{\bf\XEA\ Definition }In the rest of this paper, assume that any subcategory is full and closed under isomorphism, direct sums and direct summands.
Let $\cc^\prime$ be a subcategory of a Krull-Schmidt category $\cc$. 

(1)[I4;5.1] We call $\cc^\prime$ a {\it rejective subcategory} of $\cc$ if the inclusion functor $\cc^\prime\rightarrow\cc$ has a right adjoint $(\ )^-:\cc\rightarrow\cc^\prime$ with a counit $\epsilon^-$ such that $\epsilon^-_X$ is a monomorphism for any $X\in\cc$, and a left adjoint $(\ )^+:\cc\rightarrow\cc^\prime$ with a unit $\epsilon^+$ such that $\epsilon^+_X$ is an epimorphism for any $X\in\cc$ (compare with torsion theories \XAB).
We call $\cc^\prime$ a {\it trivial subcategory} of $\cc$ if $\cc$ is a unique rejective subcategory of $\cc$ which contains $\cc^\prime$.

(2)[I4;8.2] We call a subset ${\bf S}$ of $\ind\cc$ {\it rejectable} (resp. {\it trivial}) if ${\bf S}=\ind\cc-\ind\cc^\prime$ for some rejective (resp. trivial) subcategory $\cc^\prime$ of $\cc$.

\vskip1em{\bf\XEAA\ Example }[I4;5.3]
Let $\Lambda$ be an artin algebra and $\cc=\mod\Lambda$.

Any ring morphism $G:\Lambda\rightarrow\Gamma$ induces a faithful functor $G^*:\mod\Gamma\rightarrow\cc$, which is full if $G$ is surjective.
A subcategory $\cc^\prime$ of $\cc$ is rejective if and only if there exists a factor algebra $\Gamma$ of $\Lambda$ such that $\cc^\prime=\mod\Gamma$.
In this case, adjoint functors are given by $(\ )^+=\Gamma\otimes_\Lambda\ $ and $(\ )^-=\hom_\Lambda(\Gamma,\ )$. Thus we obtain a bijection from factor algebras of $\Lambda$ to rejective subcategories of $\cc$ defined by $\Gamma\mapsto\mod\Gamma$.

Note that it is shown in [I4;6.3] that any rejective subcategory $\cc^\prime$ of a (strict) $\tau$-category $\cc$ forms a (strict) $\tau$-category again if $\cc/[\cc^\prime]$ is artinian.

\vskip1em{\bf\XEB\ Theorem }{\it
Let $\Lambda$ be an artin algebra, $(\tt,\cc)$ a faithful torsion theory on $\mod\Lambda$ and $\cc^\prime$ a subcategory of $\cc$. Then $\cc^\prime$ is a rejective subcategory of $\cc$ if and only if there exists a morphism $G:\Lambda\rightarrow\Gamma$ of artin algebras such that $\cc^\prime=\cc\cap\quo\Gamma$, $\Gamma\in\cc^\prime$ and $\Gamma/G(\Lambda)\in\tt$, where we denote by $\quo\Gamma$ the full subcategory of $\mod\Lambda$ consisting of factor modules of $\Gamma^n$ ($n>0$). In this case, $G^*$ induces an equivalence from a faithful torsionfree class $\{X\in\mod\Gamma\ |\ G^*X\in\cc\}$ on $\mod\Gamma$ to $\cc^\prime$.}

\vskip1em{\sc Proof }
(i) We will show the `only if' part.

Let $(\ )^+:\cc\rightarrow\cc^\prime$ be a left adjoint of the inclusion functor with a unit $\epsilon^+$, and $\Gamma:=\endm_{\Lambda}(\Lambda^+)$.
Then $\epsilon^+_\Lambda\in\hom_{\Lambda}(\Lambda,\Lambda^+)$ is given by a left multiplication of an element $a\in \Lambda^+$.
Taking $\hom_\Lambda(\ ,\Lambda^+)$, we obtain a bijection $(a\cdot):\endm_\Lambda(\Lambda^+)=\Gamma\rightarrow\hom_{\Lambda}(\Lambda,\Lambda^+)=\Lambda^+$.
Thus a ring morphism $G:\Lambda\rightarrow\Gamma$ is well-defined by $xa=aG(x)$ for any $x\in\Lambda$.
Since $(a\cdot)$ is a bijection such that $(a\cdot)\circ G=(\cdot a)=\epsilon^+_\Lambda$, we can replace $\Lambda^+$ and $\epsilon^+_\Lambda$ by $\Gamma$ and $G$. Let $\Lambda\stackrel{G}{\rightarrow}\Gamma\stackrel{H}{\rightarrow}\Gamma/G(\Lambda)\rightarrow0$ be exact. Taking $\hom_\Lambda(\ ,X)$ for any $X\in\cc$, we obtain $\Gamma/G(\Lambda)\in{}^\perp\cc=\tt$.
For any $X\in\cc\cap\quo\Gamma$, take a surjection $f\in\hom_\Lambda(\Gamma^n,X)$. Since $f$ factors through the injection $\epsilon^-_X\in\hom_\Lambda(X^-,X)$, we obtain $X=X^-\in\cc^\prime$. Conversely, for any $X\in\cc^\prime$, take a surjection $f\in\hom_\Lambda(\Lambda^n,X)$. Since $f^+\in\hom_\Lambda(\Gamma^n,X)$ is surjective again, we obtain $X\in\cc\cap\quo\Gamma$. Thus $\cc^\prime=\cc\cap\quo\Gamma$ holds.

(ii) We will show the `if' part.

Fix any $X\in\cc$. Put $X^-:=\hom_\Lambda(\Gamma,X)$. Then the natural map $\epsilon^-_X\in\hom_\Lambda(X^-,X)$ is injective by $\Gamma/G(\Lambda)\in\tt$. Thus $X^-$ is a unique maximal submodule of $X$ such that $X^-\in\quo\Gamma$. Hence $X^-\in\cc\cap\quo\Gamma=\cc^\prime$ holds, and $\hom_\Lambda(\ ,X^-)\stackrel{\cdot\epsilon^-_X}{\longrightarrow}\hom_\Lambda(\ ,X)$ is an isomorphism on $\quo\Gamma\supseteq\cc^\prime$. Thus $(\ )^-:\cc\rightarrow\cc^\prime$ gives a right adjoint of the inclusion functor with a counit $\epsilon^-$.
Let $\fff:\mod\Lambda\rightarrow\cc$ be the left adjoint of the inclusion functor (\XAB), $X^+:=\fff(\Gamma\otimes_\Lambda X)\in\cc\cap\quo\Gamma=\cc^\prime$ and $\epsilon^+_X\in\hom_\Lambda(X,X^+)$ the natural map. Then $\hom_\Lambda(X^+,Y)=\hom_\Lambda(\Gamma\otimes_\Lambda X,Y)=\hom_\Lambda(X,Y^-)=\hom_\Lambda(X,Y)$ holds for any $Y\in\cc^\prime$. Thus $(\ )^+:\cc\rightarrow\cc^\prime$ gives a left adjoint of the inclusion functor with a unit $\epsilon^+$.


(iii) We will show the latter assertion.

Obviously $\ff^\prime:=\{X\in\mod\Gamma\ |\ G^*X\in\cc\}$ forms a torsionfree class on $\mod\Gamma$ and $G^*$ induces a functor $\ff^\prime\rightarrow\cc\cap\quo\Gamma=\cc^\prime$. For any $X\in\cc^\prime$, we have an isomorphism $\hom_\Lambda(\Gamma,X)=X^-\stackrel{\epsilon^-_X}{\rightarrow}X$ by (ii). Since $\Gamma$ is a $(\Lambda,\Lambda)$-bimodule, we can regard $X^-=\hom_\Lambda(\Gamma,X)$ as a $\Gamma$-modules such that $G^*(X^-)=X$. Thus $G^*:\ff^\prime\rightarrow\cc^\prime$ is dense.
Finally we will show that $G^*$ is full faithful on $\ff^\prime$. For any $X,Y\in\ff^\prime$, take an exact sequence $\Gamma^m\rightarrow\Gamma^n\rightarrow X\rightarrow0$ in $\mod\Gamma$. Taking $\hom_{\Lambda}(\ ,Y)$ and $\hom_{\Gamma}(\ ,Y)$, we obtain $\hom_\Lambda(X,Y)=\hom_\Gamma(X,Y)$ by $Y=\hom_\Lambda(\Gamma,Y)$.\rule{5pt}{10pt}

\vskip1em{\bf\XEC\ }
Let $\cc$ be a $\tau$-category, $\cc^\prime$ a subcategory of $\cc$, and $\overline{\cc}:=\cc/[\cc^\prime]$ the factor category.
Then we can regard $\ind\cc$ as a disjoint union of $\ind\overline{\cc}$ and $\ind\cc^\prime$ naturally.
By [I4;1.4], $\overline{\cc}$ forms a $\tau$-category again, and $\AR(\overline{\cc})$ is obtained by deleting vertices in $\ind\cc^\prime$ from $\AR(\cc)$. Thus we easily obtain the terms $\theta^\pm_{\overline{\cc}}$ and $\tau^\pm_{\overline{\cc}}$ of left and right $\tau$-sequences in $\overline{\cc}$.
Notice that non-trivial subcategory $\cc^\prime$ of $\cc$ is rejective if any subcategory of $\cc$ containing $\cc^\prime$ is trivial except $\cc^\prime$. Thus we can use both (1) and (2) below to check the rejectivity.

\vskip1em{\bf Theorem }{\it
Let $\cc^\prime$ be a rejective subcategory of a strict $\tau$-category $\cc$. Assume that $\overline{\cc}:=\cc/[\cc^\prime]$ is artinian.

(1)[I4;8.2.1] $\cc^\prime$ is a trivial subcategory of $\cc$ if and only if $\overline{\cc}(B,A)=0$ holds for any $A\in\ii\cc$ and $B\in\pp\cc$ if and only if the condition below is satisfied.

For any $X\in\ind\overline{\cc}\cap\ii\cc$, put $Y_0:=X$, $Y_1:=\theta^+_{\overline{\cc}}X$ and $Y_i:=(\theta^+_{\overline{\cc}}Y_{i-1}-\tau^+_{\overline{\cc}}Y_{i-2})_+$ for $i\ge2$. Then $Y_i\in\zzz(\ind\overline{\cc}-\pp\cc)$ holds for any $i\ge0$.

(2)[I4;8.2.2] $\cc^\prime$ is a rejective subcategory of $\cc$ if and only if $\mu^-_A$ is a monomorphism and $\mu^+_B$ is an epimorphism in $\overline{\cc}$ for any $A\in\ind\overline{\cc}-\ii\cc$ and $B\in\ind\overline{\cc}-\pp\cc$ if and only if (i) and (ii) below are satisfied.


\strut\kern1em(i) For any $X\in\ind\overline{\cc}-\ii\cc$, put
$Y_0:=\theta^-_{\overline{\cc}}X$, $Y_1:=\theta^+_{\overline{\cc}}
\theta^-_{\overline{\cc}}X-X$ and $Y_i:=\theta^+_{\overline{\cc}}Y_{i-1}-
\tau^+_{\overline{\cc}}Y_{i-2}$ for $i\ge2$.
Then $Y_i\in\nnn\ind\overline{\cc}$ holds for any $i\ge0$.

\strut\kern1em(ii) For any $X\in\ind\overline{\cc}-\pp\cc$, put
$Y_0:=\theta^+_{\overline{\cc}}X$, $Y_1:=\theta^-_{\overline{\cc}}
\theta^+_{\overline{\cc}}X-X$ and $Y_i:=\theta^-_{\overline{\cc}}Y_{i-1}-
\tau^-_{\overline{\cc}}Y_{i-2}$ for $i\ge2$.
Then $Y_i\in\nnn\ind\overline{\cc}$ holds for any $i\ge0$.}

\vskip1em{\bf\XECA\ Corollary (DK Rejection Lemma) }{\it
Let $\cc^\prime$ be a rejective subcategory of a strict $\tau$-category $\cc$. Assume that $\ind\cc-\ind\cc^\prime=\{ X\}$ and $\cc/[\cc^\prime]$ is artinian. Then $\cc^\prime$ is a rejective subcategory of $\cc$ if and only if $X\in\pp\cc\cap\ii\cc$ holds.}

\vskip1em{\bf\XED\ Example }
(1) By \XECA, singleton sets $\{4\}$, $\{6\}$, $\{8\}$, $\{24\}$ in \XDE(1) and $\{1\}$,$\{4\}$ in \XDE(3) are rejectable.

(2) We can easily check that
{\scriptsize$\begin{array}{rcl}11\rightarrow8\rightarrow6&\rightarrow&3\\
\downarrow&&\downarrow\\
4&\rightarrow&24\rightarrow22\rightarrow19\end{array}$}
in \XDE(2) is rejectable by \XEC(2) (or (1)).
Moreover, the AR quiver $\AR(\cc^\prime)$ of the corresponding rejective subcategory $\cc^\prime$ is the following:
{\scriptsize\[\begin{array}{cccccccccccccccccc}
&&&&&&&&\Q{7}&&&&\Q{1}&&&\\
&&&&&&&\QUR&&\QDR&&\QUR&&&&&\\
&&\Q{20}&&&&\Q{14}&&&&\Q{5}&&&&&&\\
&\QUR&&\QDR&&\QUR&&\QDR&&\QUR&&\QDR&&&&&\\
\Q{2}&\QR&\Q{21}&\QR&\Q{18}&\QR&\Q{15}&\QR&\Q{12}&\QR&\Q{9}&\QR&\Q{2}&\ \ \ \ &\Q{13}&\QR&\Q{17}\\
&\QDR&&\QUR&&\QDR&&\QUR&&\QDR&&\QUR&&&&&\\
&&\Q{25}&&&&\Q{16}&&&&\Q{10}&&&&&&\\
&\QUR&&\QDR&&\QUR&&&&&&&&&&&\\
\Q{26}&&&&\Q{23}&&&&&&&&&&&&
\end{array}\ \ \ \ \ \begin{array}{c}
\tau^+(20)=10\\
Q^p=\{7,13,17,25,26\}\\
Q^i=\{1,5,13,17,23\}\end{array}\]}

(3) We can easily check that {\scriptsize$20\rightarrow17\rightarrow12\rightarrow8$} and {\scriptsize$\begin{array}{lll}21\rightarrow&17\rightarrow&12\\
\downarrow&\downarrow&\downarrow\\
18\rightarrow&14\rightarrow&10\rightarrow6\end{array}$}
and so on in \XDE(4) are rejectable

\vskip1em{\bf\XEE\ }
Let $\cc^\prime$ be a rejective subcategory of a strict $\tau$-category $\cc$. If $\overline{\cc}:=\cc/[\cc^\prime]$ is artinian, then $\overline{\cc}$ forms a strict $\tau$-category again by [I4;6.1] (cf. \XEC(2)). In particular, we obtain the following result by \XBA.

\vskip1em{\bf Proposition }{\it
Let $\cc$ be a faithful torsionfree class over an artin algebra $\Lambda$, $\cc^\prime$ a rejective subcategory of $\cc$ and $\overline{\cc}:=\cc/[\cc^\prime]$. If $\#\ind\overline{\cc}<\infty$, then $\overline{\cc}$ is equivalent to a faithful torsionfree class over some artin algebra $\Lambda^\prime$.}

\vskip1em{\bf\XEEA\ }The above result \XEE\ holds even if we drop the assumption $\#\ind\overline{\cc}<\infty$. See the author's forthcoming papers.

\vskip1em{\footnotesize
\begin{center}
{\bf References}
\end{center}

[A] M. Auslander: Representation theory of Artin algebras. I, II. Comm. Algebra 1 (1974), 177--268; ibid. 1 (1974), 269--310.

[AB] M. Auslander, M. Bridger: Stable module theory. Memoirs of the American Mathematical Society, No. 94 American Mathematical Society, Providence, R.I. 1969 146 pp.

[AR1] M. Auslander, I. Reiten: Representation theory of Artin algebras. III. Almost split sequences. Comm. Algebra 3 (1975), 239--294.

[AR2] M. Auslander, I. Reiten: $k$-Gorenstein algebras and syzygy modules. J. Pure Appl. Algebra 92 (1994), no. 1, 1--27.

[AR3] M. Auslander, I. Reiten: Syzygy modules for Noetherian rings. J. Algebra 183 (1996), no. 1, 167--185.

[ARS] M. Auslander, I. Reiten, S. O. Smal:
Representation theory of Artin algebras. 
Cambridge Studies in Advanced Mathematics, 36. 
Cambridge University Press, Cambridge, 1995. 

[AS] M. Auslander, S. O. Smalo: Almost split sequences in subcategories. J. Algebra 69 (1981), no. 2, 426--454. 

[As] I. Assem: Tilting theory---an introduction. Topics in algebra, Part 1 (Warsaw, 1988), 127--180, Banach Center Publ., 26, Part 1, PWN, Warsaw, 1990.

[B] H. Bass: On the ubiquity of Gorenstein rings. Math. Z. 82 1963 8--28. 

[Bj] J-E. Bjork: The Auslander condition on Noetherian rings. Seminaire d'Algebre Paul Dubreil et Marie-Paul Malliavin, 39eme Annee (Paris, 1987/1988), 137--173, Lecture Notes in Math., 1404, Springer, Berlin, 1989. 

[BG] K. Bongartz, P. Gabriel:
Covering spaces in representation-theory. 
Invent. Math. 65 (1981/82), no. 3, 331--378.

[Br] S. Brenner: A combinatorial characterization of finite Auslander-Reiten quivers, (Ottawa, Ont., 1984), 13--49, Lecture Notes in Math., 1177, Springer, Berlin-New York, 1986.

[C] J. Clark: Auslander-Gorenstein rings for beginners. International Symposium on Ring Theory (Kyongju, 1999), 95--115, Trends Math., Birkhauser Boston, Boston, MA, 2001.

[CR] C. W. Curtis, I. Reiner: Methods of representation theory. Vol. I. With applications to
finite groups and orders. Reprint of the 1981 original. Wiley Classics Library. A Wiley-Interscience Publication. John
Wiley \& Sons, Inc., New York, 1990.

[DK] Y. A. Drozd, V. V. Kiri\v cenko: The quasi-Bass orders. (Russian) Izv. Akad. Nauk SSSR Ser.
Mat. 36 (1972), 328--370.

[DKR] Y. A. Drozd, V. V. Kiri\v cenko, A. V. Ro\u\i ter: Hereditary and Bass orders. (Russian) Izv. Akad. Nauk SSSR Ser. Mat. 31 1967 1415--1436.

[FGR] R. M. Fossum, P. Griffith, I. Reiten: Trivial extensions of abelian categories. Homological algebra of trivial extensions of
abelian categories with applications to ring theory. Lecture Notes in Mathematics, Vol. 456. Springer-Verlag, Berlin-New York, 1975.

[GM] E. L. Green, R. Martinez Villa: Koszul and Yoneda algebras. Representation theory of algebras (Cocoyoc, 1994), 247--297, CMS Conf. Proc., 18, Amer. Math. Soc., Providence, RI, 1996.

[Ha] D. Happel: Triangulated categories in the representation theory of finite-dimensional algebras. 
London Mathematical Society Lecture Note Series, 119. 
Cambridge University Press, Cambridge, 1988.

[H1] M. Hoshino: Tilting modules and torsion theories. Bull. London Math. Soc. 14 (1982), no. 4, 334--336. 

[H2] M. Hoshino: On dominant dimension of Noetherian rings. Osaka J. Math. 26 (1989), no. 2, 275--280.

[HPR] D. Happel, U. Preiser, C. M. Ringel: Vinberg's characterization of Dynkin diagrams using subadditive functions with application to $D{\rm Tr}$-periodic modules. Representation theory, II (Proc. Second Internat. Conf., Carleton Univ., Ottawa, Ont., 1979), pp. 280--294, Lecture Notes in Math., 832, Springer, Berlin, 1980.


[HN] H. Hijikata, K. Nishida:
Bass orders in nonsemisimple algebras. J. Math. Kyoto Univ. 34 (1994),
no. 4, 797--837.

[I1] O. Iyama: A generalization of rejection lemma of Drozd-Kirichenko. J. Math. Soc. Japan 50 (1998), no. 3, 697--718.

[I2] O. Iyama: Some categories of lattices associated to a central idempotent. J. Math. Kyoto Univ. 38 (1998), no. 3, 487--501.

[I3] O. Iyama: $\tau$-categories I: Ladders, to appear in Algebras and Representation Theory.

[I4] O. Iyama: $\tau$-categories II: Nakayama pairs and rejective subcategories, to appear in Algebras and Representation Theory.

[I5] O. Iyama: $\tau$-categories III: Auslander orders and Auslander-Reiten quivers, to appear in Algebras and Representation Theory.

[I6] O. Iyama: Representation theory of orders. "Algebra - Representation Theory" (Edited by K. W. Roggenkamp, M. Stefanescu), Nato-Science-Series, Kluwer Academic Publishers, 63-96 (2001)

[I7] O. Iyama: A proof of Solomon's second conjecture on local zeta functions of orders, J. Algebra 259 (2003), no. 1, 119-126.

[I8] O. Iyama: Finiteness of Representation dimension, Proc. Amer. Math. Soc. 131 (2003), no. 4, 1011-1014.

[I9] O. Iyama: Symmetry and Duality on $n$-Gorenstein rings, to appear in Journal of Algebra.

[IT1] K. Igusa, G. Todorov:
Radical layers of representable functors. J. Algebra 89 (1984),
no. 1, 105--147.

[IT2] K. Igusa, G. Todorov:
A characterization of finite Auslander-Reiten quivers. J.
Algebra 89 (1984), no. 1, 148--177.

[IT3] K. Igusa, G. Todorov: A numerical characterization of finite Auslander-Reiten quivers. Representation theory, I (Ottawa, Ont., 1984), 181--198, Lecture Notes in Math., 1177, Springer, Berlin, 1986.

[Ro] K. W. Roggenkamp:
Lattices over orders. II. Lecture Notes in Mathematics, Vol. 142
Springer-Verlag, Berlin-New York 1970.

[R1] W. Rump: Almost abelian categories. Cahiers Topologie Geom. Differentielle Categ. 42 (2001), no. 3, 163--225.

[R2] W. Rump: $*$-modules, tilting, and almost abelian categories. Comm. Algebra 29 (2001), no. 8, 3293--3325.

[R3] W. Rump: Ladder functors with an application to representation-finite Artinian rings. An. \c Stiin\c t. Univ. Ovidius Constan\c ta Ser. Mat. 9 (2001), no. 1, 107--123.

[R4] W. Rump: Lattice-finite rings, preprint.

[R5] W. Rump: The category of lattices over a lattice-finite ring, preprint.

[RV] I. Reiten, M. Van den Bergh:
Two-dimensional tame and maximal orders of finite representation type. 
Mem. Amer. Math. Soc. 80 (1989).

[RVo] C. M. Ringel, D. Vossieck, Hammocks,
Proc. London Math. Soc. (3) 54 (1987), 216--246.

[S] D. Simson: Linear representations of partially ordered sets and vector space categories. Algebra, Logic and Applications, 4. Gordon and Breach Science Publishers, Montreux, 1992.

[T] H. Tachikawa: Quasi-Frobenius rings and generalizations. Lecture Notes in Mathematics, Vol. 351. Springer-Verlag, Berlin-New York, 1973.
}

\vskip.5em{\footnotesize
{\sc Department of Mathematics, Himeji Institute of Technology, Himeji, 671-2201, Japan}

{\it iyama@sci.himeji-tech.ac.jp}}

\end{document}